\input amstex
\documentstyle{amsppt}
\document
\magnification=1200
\NoBlackBoxes
\vsize17cm
\nologo
\def\P{\bold{P}}
\def\Q{\bold{Q}}
\def\R{\bold{R}}
\def\C{\bold{C}}
\def\Z{\bold{Z}}
\def\H{\bold{H}}
\def\cT{{\Cal{T}}}
\def\cC{{\Cal{C}}}
\def\cA{{\Cal{A}}}
\def\fA{\bold{A}}
\def\cI{{\Cal{I}}}
\def\cJ{{\Cal{J}}}


\bigskip

\centerline{\bf MODULAR SHADOWS}

\smallskip

\centerline{\bf AND THE L\'EVY--MELLIN $\infty$--ADIC TRANSFORM}

\medskip

\centerline{\bf Yuri I.~Manin, Matilde Marcolli}

\medskip

\centerline{\it Max--Planck--Institut f\"ur Mathematik, Bonn, Germany,}

\centerline{\it and Northwestern University, Evanston, USA}

\bigskip

{\bf Abstract.} This paper continues the study of the
structures induced on the
``invisible boundary'' of the modular tower and extends some
results of [MaMar1].
We start with a systematic formalism
of pseudo--measures generalizing the well--known theory
of modular symbols for $SL(2)$. These pseudo--measures, and
the related integral formula which we call the L\'evy--Mellin
transform, can be considered
as an ``$\infty$--adic'' version of Mazur's $p$--adic measures
that have been introduced in the seventies in the theory of
$p$--adic interpolation of the Mellin transforms of cusp forms,
cf. [Ma2]. A formalism of iterated  L\'evy--Mellin
transform in the style of [Ma3] is sketched.
Finally, we discuss the invisible boundary from the
perspective of non--commutative geometry.
\bigskip

\centerline{\bf 0. Introduction}

\medskip

When the theory of modular symbols  for the $SL(2)$--case had been conceived
in the 70's (cf. [Ma1], [Ma2], [Sh1], [Sh2]),
it was clear from the outset that it dealt with the
Betti homology of some basic moduli spaces (modular curves,
Kuga varieties, $\overline{M}_{1,n}$, and alike),
whereas the theory of modular forms involved the
de Rham and Hodge cohomology of the same spaces.

\smallskip

However, the combinatorial skeleton of the formalism of modular
symbols is so robust, depending essentially only on the
properties of continued fractions, that other
interpretations and connections naturally suggest themselves.

\smallskip

In this paper, we develop the approach to the
modular symbols which treats them as a special case 
of some structures supported by the
``invisible boundary'' of the tower of classical modular curves
along the lines of [MaMar1], [MaMar2].

\smallskip

Naively speaking, this boundary is (the tower of) quotient space(s) of
$\P^1(\R )$ modulo (finite index subgroups of) $PSL(2,\Z )$.
The part of it consisting of orbits of ``cusps'',  $\P^1(\Q )$,
has a nice algebraic geometric description, but
irrational points are not considered in
algebraic geometry, in particular, since the action of $PSL(2,\Z )$
is highly non--discrete. This is why we call this part
``invisible''.

\smallskip 

``Bad quotients'' of this type can be efficiently studied 
using tools of Connes'
non--commutative geometry.
Accordingly, the Theorem 4.4.1 of [MaMar1] identified
the modular homology complex with (a part of)
Pimsner's exact sequence in the $K_*$--theory of the reduced crossed
product algebra $C(\P^1(\R ))\rtimes PSL(2,\Z )$.
Moreover, the Theorem 0.2.2 of the same 
paper demonstrated the existence 
of a version of Mellin transform (from modular forms
to Dirichlet series) where the integrand was supported
by the real axis rather than upper complex half--plane.

\smallskip

In [MaMar2], these and similar results were put in connection
with the so called ``holography'' principle in modern theoretical physics.
According to this principle, quantum field theory on a space may
be faithfully reflected by an appropriate theory on
the boundary of this space. When this boundary,
rather than the interior, is interpreted as our observable space--time,
one can proclaim that the ancient Plato's cave metaphor 
is resuscitated in this sophisticated guise.
This metaphor motivated the title of the present paper.  

\smallskip

Here is a review of its contents.

\smallskip

The sections 1--4 address the first of the two basic themes:

\smallskip

(i) {\it How does the holomorphic geometry of the upper complex
half plane project itself onto the invisible boundary?}

\smallskip

In the first and the second sections, we introduce and develop a 
formalism of general pseudo--measures. They can be defined as
finitely additive functions with values in an abelian group $W$
supported by the 
Boolean algebra generated by segments with rational
ends in $\P^1(\R )$. Although the definitions (and proofs) 
are very elementary,
they capture some essential properties of the modular boundary.

\smallskip
 
In particular, in subsections 1.9--1.11
we show that the {\it generalized Dedekind symbols}
studied by Sh.~Fukuhara in [Fu1], [Fu2] are essentially
certain sequences of pseudo--measures. Example 1.12
demonstrates that modular symbols are pseudo--measures.
(In fact, much more general integrals along geodesics connecting cusps
produce pseudo--measures; this is why we use the imagery
of ``projecting  the holomorphic geometry of the upper complex
half plane onto the invisible boundary''.)
Finally, subsection 2.3.1 establishes that {\it rational period
functions} in the sense of [Kn1], [Kn2], [A], [ChZ],
are pseudo--measures as well.

\smallskip

In the third section, we define the L\'evy--Mellin
transform of a pseudo--measure and prove the Theorem 3.4
generalizing the Theorem 0.2.2 of [MaMar1].
This shows that some specific L\'evy functions
involving pseudo--measures can serve as an efficient
replacement of (non--existent) restrictions
of cusp forms to the real boundary of $H$.

\smallskip

In the fourth section, a non--commutative 
version of pseudo--measures is developed. The new formalism
was suggested by the theory of iterated integrals
of modular forms introduced in [Ma3], [Ma4].
The iterated L\'evy--Mellin transform appears naturally in this context
producing some interesting multiple Dirichlet series
related to but different from those discussed in [Ma3].
We hope to return to this subject later.

\smallskip

The second theme developed in the fifth section is:

\smallskip

(ii) {\it With what natural structure(s) is the invisible
boundary endowed?}

\smallskip

We start with a reformulation of pseudo--measures in terms of currents 
on the tree $\cT$ of $PSL(2,\Z)$.  The boundary action of
$PSL(2,\Z)$ can be best visualized in terms of the space of ends of 
this tree. This set of endpoints is a compact Hausdorff
space which maps continuously to $\P^1(\R)$ through a natural map that 
is $1:1$ on the irrational points and $2:1$ on the rationals. 

\smallskip

In terms of noncommutative geometry, this boundary action is described 
by a crossed product $C^*$-algebra, and $W$--valued pseudo--measures can be
interpreted as homomorphisms from the $K_0$ of the crossed product
$C^*$-algebra to $W$.

\smallskip 

We also introduce another $C^*$-algebra naturally associated to the
boundary action, described in terms of the generalized Gauss shift
of the continued fraction expansion. We show that this can be 
realized as a subalgebra of the crossed product algebra of
the action of $PSL(2,\Z)$ on $\partial\cT$.

\smallskip

We then consider an extension of pseudo--measures to limiting
pseudo--measures that parallels the notion of 
limiting modular symbols introduced in [MaMar1].
The limiting modular symbols can then be realized as limits 
of pseudo--measures associated to the ordinary modular symbols
or as averages of currents on the tree of $PSL(2,\Z)$.

\bigskip  

\centerline{\bf 1. Pseudo--measures: commutative case}

\medskip

In the following
we collect the basic facts from the theory of Farey series
in the form that we will use throughout the paper.

\medskip

{\bf 1.1. Conventions.} We will consider $\Q\subset \R\subset \C$
as points of an affine line with a fixed coordinate, say, $z$.
Completing this line by one point $\infty$, we get
points of the projective line $\P^1(\Q )\subset \P^1(\R )\subset \P^1(\C )$.
{\it Segments} of $\P^1(\R )$ are defined as non--empty 
connected subsets of $\P^1(\R )$. The boundary of each segment generally
consists of an (unordered) pair of points in $\P^1(\R )$. In marginal
cases, the boundary might be empty or consist of one point
in which case we may speak of an {\it improper segment}. A proper segment
is called {\it rational} if its ends are in $\P^1(\Q )$.
It is called infinite if $\infty$ is in its closure,
otherwise it is called finite.

\medskip

{\bf 1.2. Definition.} {\it A pseudo--measure on
${\P}^1(\R )$ with values in a commutative group  (written additively)
$W$ is a function $\mu :\, {\P}^1(\Q )\times {\P}^1(\Q )\to W$
satisfying the following conditions: for any $\alpha ,\beta ,\gamma \in
{\P}^1(\Q )$},
$$
\mu (\alpha ,\alpha )=0,\quad \mu (\alpha ,\beta ) +\mu (\beta ,\alpha )=0\,,
\eqno(1.1)
$$
$$
\mu (\alpha ,\beta ) + \mu (\beta ,\gamma ) + \mu (\gamma ,\alpha )=0\, .
\eqno(1.2)
$$

\medskip

There are two somewhat different ways to look at $\mu$
as a version of measure.

\smallskip

(i) We can uniquely extend the map
$(\alpha ,\beta )\mapsto \mu (\alpha ,\beta )$ to a finitely additive function on the
Boolean algebra consisting of finite unions of {\it positively
oriented from $\alpha$ to $\beta$} rational segments. On improper
segments, in particular $points$ and the whole $\P^1 (\R)$,
this function vanishes.

\smallskip

Here positive orientation is defined by the increasing $z$. Thus, the ordered
pair $(1,2)$ corresponds to the segment $1\le z \le 2$
(one or both ends can be excised, the pseudo--measure remains the same).
However the pair $(2,1)$ corresponds  to the union 
$\{2\le z \le \infty \}\cup \{-\infty \le z \le 1\}$ in the traditional notation
(which can be somewhat misleading since $-\infty = \infty$ in
$\P^1(\R ).$) Thus, $(2,1)$ is an infinite segment.

\smallskip

(ii) Another option is to restrict oneself to finite segments
and assign to $(\alpha ,\beta )$
the segment with these ends {\it oriented from $\alpha$ to $\beta$.}

\smallskip 

We will freely use both viewpoints, and allow ourselves some
laxity about which end belongs to our
segment and which not whenever it is not essential. 

\smallskip

Moreover, working with pseudo--measures as purely combinatorial 
objects, we will simply identify {\it segments} and {\it ordered pairs} 
$(\alpha ,\beta )
\in \P^1 (\Q )$. We will call $\alpha$ (resp. $\beta$)
the {\it ingoing} (resp. {\it outgoing}) end. 

\medskip

{\bf 1.3. The group of pseudo--measures.} $W$--valued pseudo--measures form a commutative group $M_W$ 
with the composition law 
$$
(\mu_1 +\mu_2)(\alpha ,\beta ):=
\mu_1 (\alpha ,\beta ) + \mu_2(\alpha ,\beta ).
\eqno(1.3)
$$
If $W$ is an $A$--module over a ring $A$, the pseudo--measures
form an $A$--module as well.

\medskip

{\bf 1.4. The universal pseudo--measure.} Let $\Z [\P^1(\Q )]$ be a free abelian group freely generated
by $\P^1(\Q )$, and $\nu :\, \P^1(\Q )\to \Z[\P^1(\Q )]$
the tautological map. Put $\mu^U (\alpha ,\beta ) :=
\nu (\beta )-\nu (\alpha ).$ Clearly, this is a pseudo--measure 
taking its values in the subgroup  $\Z [\P^1(\Q )]_0$,
kernel of the augmentation map $\sum_i m_i\nu (\alpha_i)\mapsto
\sum_i m_i$. This pseudo--measure is universal in the following sense: 
for any pseudo--measure  $\mu$ with values in $W$,
there is a homomorphism $w: \,\Z [\P^1(\Q )]_0 \to W$
such that $\mu = w\circ \mu^U$. In fact, we have to put
$w(\nu (\beta )-\nu (\alpha ))=\mu (\alpha ,\beta ).$

\medskip

{\bf 1.5. The action of $GL(2,\Q )$.}  The group $GL(2,\Q )$ acts 
from the left upon $\P^1(\R)$
by fractional linear automorphisms $z\mapsto gz$, mapping
$\P^1(\Q )$ to itself. It acts on pseudo--measures
with values in $W$ from the right by the formula
$$
(\mu g)(\alpha ,\beta ):= \mu (g(\alpha ) ,g(\beta )) .
\eqno(1.4)
$$
This action is compatible with the structures described in 1.3.

\smallskip

The map $z\mapsto -z$ defines an involution on $M_W$,
whose invariant, resp., antiinvariant points can be called
even, resp. odd measures.

\medskip

{\bf 1.6. Primitive segments and primitive chains.} 
A segment $I$ is called {\it primitive}, 
if it is rational,
and if its ends are of the form $\left(\dfrac{a}{c},\dfrac{b}{d}\right)$,
$a,b,c,d \in \Z$ such that $ad-bc=\pm 1$. In other words,
$I=(g(\infty ), g(0))$ where  $g\in GL(2,\Z )$ and the group
$GL(2)$ acts upon $\P^1$ by fractional
linear transformations. 

\smallskip

If $\roman{det}\,g=-1$,
we can simultaneously change signs of the entries of  the second
column. The segment $I=(g(\infty ), g(0))$ will remain the same.
The intersection of the stationary subgroups of $\infty$
and $0$ in $SL(2,\Z )$ is $\pm id$. Hence
the set of oriented primitive segments is a principal
homogeneous space over $PSL(2,\Z )$.

\smallskip

The rational ends of $I$ above 
are written in lowest
terms, whereas $\infty$ must be written as $\dfrac{\pm 1}{0}$. Signs
of the numerators/denominators are generally not normalized
(can be inverted simultaneously),
but it is natural to use $\dfrac{1}{0}$ for $+\infty$ and
$\dfrac{-1}{0}$ for $-\infty$ whenever we imagine our pairs 
as ends of oriented segments.  

\smallskip

Primitive segments with one infinite end are thus
$(-\infty ,m)$ and $(n,\infty )$, $m,n\in \Z$.
A segment with finite ends $(\alpha ,\beta )$ is primitive iff
$|\alpha - \beta |=n^{-1}$ for some $n\in \Z$, $n\ge 1$.

\smallskip

Shifting a primitive segment $I=(\alpha ,\beta )$ by any integer,
or changing signs to $-I =(-\alpha ,-\beta )$  we get a primitive
segment. Hence any finite primitive segment
can be shifted into $[0,1]$, and any primitive segment $I$ of length
(Lebesgue measure) $|I|\le \dfrac{1}{2}$ after an appropriate shift lands entirely either in $\left[0,\dfrac{1}{2}\right]$,
or in $\left[\dfrac{1}{2}, 1\right]$. Moreover, if $I=(\alpha ,\beta )\subset
\left[0, \dfrac{1}{2}\right]$ is primitive, then 
  $1-I=:(1-\alpha ,1-\beta )\subset \left[\dfrac{1}{2}, 1\right]$ is primitive,
and vice versa. 

\smallskip

Generally, let $\alpha ,\beta \in \P^1(\Q )$. Let us call
{\it a primitive chain of length $n$ connecting $\alpha$ to $\beta$}
any non--empty fully ordered family of proper primitive segments
$I_1,\dots ,I_n$ such that the ingoing end of $I_1$
is $\alpha$, outgoing end of $I_n$ is $\beta$, 
and for each $1\le k\le n-1$, the outgoing end of
$I_k$ coincides with the ingoing end of $I_{k+1}$.
The numeration of the segments is not a part of the
structure, only their order is. We will call
chains with $\alpha =\beta$ {\it primitive loops},
and allow one improper segment $(\alpha ,\alpha )$ 
to be considered as a primitive loop of length $0$. 

\smallskip

This  notion is covariant with respect to the $SL(2,\Z )$--action:
if $I_1,\dots ,I_n$ is a primitive chain connecting $\alpha$ to $\beta$,
then for any $g\in SL(2,\Z )$,
$g(I_1),\dots ,g(I_n)$ is a primitive chain connecting $g(\alpha )$ to 
$g(\beta )$.

\medskip

{\bf 1.6.1. Lemma.} {\it (a) If $I_1$, $I_2$ are two open primitive 
segments and at least one of them is finite,
then either $I_1\cap I_2=\emptyset$, or one of them is contained
in another.

\smallskip

(b) Any two points $\alpha ,\beta \in \P^1(\Q )$ can be connected
by a primitive chain.}

\smallskip

{\bf Proof.} This is well known. We reproduce an old argument
showing (b) from [Ma1] in order to fix some notation.
Consider the following sequence of {\it normalized convergents to $\alpha$} 
$$
\frac{p_{-1}(\alpha )}{q_{-1}(\alpha )}:=\frac{1}{0}=\infty ,\ 
\frac{p_0(\alpha )}{q_0(\alpha )}:=
\frac{p_0}{1},\ \dots \ ,\frac{p_n(\alpha )}{q_n(\alpha )} =\alpha
\eqno(1.5)
$$
Here $\alpha\in \Q$, $p_0:=[\alpha ]$ the integer part of $\alpha$,
and convergents are calculated from
$$
\alpha = k_0(\alpha )+[k_1(\alpha ),\dots ,k_n(\alpha )]:= p_0 +
\frac{1}{k_1(\alpha )+\frac{1}{k_2(\alpha )+\dots \frac{1}{k_n(\alpha )}}}
$$
with $1\le k_i(\alpha )\in \Z$ for $i\ge 1$
and $k_n(\alpha )\ge 2$ whenever
$\alpha \notin \Z$ so that $n\ge 1$.
 
\smallskip

 The sequence
$$
I_k(\alpha ):= \left(\frac{p_k(\alpha )}{q_k(\alpha )},
\frac{p_{k+1}(\alpha )}{q_{k+1}(\alpha )}\right)
\eqno(1.6)
$$
is a primitive chain connecting $\infty$ to $\alpha$.

\smallskip

Applying this construction to $\beta$ and reversing the sequence of
ends (1.5), we get a primitive chain connecting $\beta$ to $\infty$.
Joining chains from $\alpha$ to $\infty$ and from
$\infty$ to $\beta$, we get a chain from $\alpha$ to $\beta$.

\smallskip

We put also
$$
g_k(\alpha ):=\left(\matrix p_k(\alpha ) & (-1)^{k+1}p_{k+1}(\alpha )\\
q_k(\alpha )& (-1)^{k+1}q_{k+1}(\alpha ) 
\endmatrix\right) \in SL (2,\Z)\, 
\eqno(1.7)
$$
so that
$$
I_k(\alpha ) = (g_k(\alpha )(\infty ), g_k(\alpha )(0))\,.
\eqno(1.8)
$$

\medskip

{\bf 1.7. Corollary.} {\it Each pseudo--measure is 
completely determined by its values on primitive 
segments.}

\smallskip

In fact, for any $\alpha$, $\beta$
and any primitive chain $I_1,\dots ,I_n$
connecting $\alpha$ to $\beta$, we must have
$$
\mu (\alpha ,\beta )= \sum_{k=1}^n \mu (I_k).
\eqno(1.9)
$$ 
 
\smallskip

{\bf 1.7.1. Definition.}  {\it A pre--measure $\widetilde{\mu}$ is
any $W$--valued function  defined on
 primitive segments and
satisfying relations (1.1) and (1.2)  for all
primitive chains of length $\le$ 3.}

\smallskip

Clearly, restricting a pseudo--measure to primitive segments, we
get  a pre--measure. We will prove the converse statement. 

\medskip

{\bf 1.8. Theorem.} {\it Each pre--measure $\widetilde{\mu}$
can be uniquely extended to a pseudo--measure $\mu$.}

\smallskip

{\bf Proof.}  If such a pseudo--measure exists, it is defined
by the (family of) formula(s) (1.9):
$$
\mu (\alpha ,\beta ):= \sum_{k=1}^n \widetilde{\mu} (I_k).
\eqno(1.10)
$$ 
We must only check that this prescription is well defined, that is,
does not depend on the choice of  $\{I_k\}$. Our argument below
is somewhat more elaborate than what is strictly needed
here. Its advantage is that it can be applied without changes
to the proof of the Theorem 4.3 below, which is a non--commutative
version of the Theorem 1.8.

\smallskip
 
We will first of all define four types of
{\it elementary moves}
which transform a primitive chain $\{I_k\}$ connecting $\alpha$ to $\beta$
to another such primitive chain  
without changing the  r.h.s. of (1.10).

\smallskip

(i) Choose $k$ (if it exists) such that $I_k=(\gamma ,\delta )$,
$I_{k+1}=(\delta ,\gamma )$, and delete $I_k,I_{k+1}$ from the chain.

\smallskip

This does not change (1.10) because 
$\widetilde{\mu} (\gamma ,\delta )+\widetilde{\mu}(\delta ,\gamma )=0$.

\smallskip

(ii) A reverse move: choose a point $\gamma$ which is
the outgoing end of some $I_k$ (resp. $\gamma =\alpha$),
choose a primitive segment
$(\gamma ,\delta )$, and insert the pair $(\gamma ,\delta ),
(\delta ,\gamma )$ right after $I_k$ (resp. before $I_1$.) 

\smallskip

This move can also be
applied to the empty loop connecting $\gamma$ to $\gamma$, then
it produces a chain of length two. Again, application of such a move
is compatible with (1.10).

\smallskip

(iii) Choose $I_k,I_{k+1},I_{k+2}$ (if they exist) such that
these three segments have the form $(\gamma_1,\gamma_2),
(\gamma_2,\gamma_3),(\gamma_3,\gamma_1)$, and delete
them from the chain.

\smallskip

This move is compatible with (1.10) as well because
we have postulated (1.2) on such chains.

\smallskip

(iv) A reverse move:  choose a point $\gamma$ which is
the outgoing end of some $I_k$ (resp. $\gamma =\alpha$), and
insert any triple of segments as above  
right after $I_k$ (resp. before $I_1$.)

\smallskip

Now we will show that 

\smallskip

(*) {\it any primitive loop, that is, a primitive 
chain with $\alpha =\beta$,
can be transformed into an empty loop by a sequence
of elementary moves.} 

\smallskip

Suppose that we know this.
If $I_1,\dots ,I_n$,  $J_1,\dots ,J_m$ are two chains
connecting  $\alpha$ to $\beta$, we can produce from them
a loop connecting $\alpha$ to $\alpha$ by putting
after $I_n$ the segments $J_m,\dots ,J_1$ with reversed
orientations. The r.h.s. part of (1.10) calculated for
this loop must be zero because it is zero for the empty
loop. Hence $I_1,\dots ,I_n$ and $J_1,\dots ,J_m$
furnish the same r.h.s. of (1.10).

\smallskip

We will establish (*)
by induction on the length of a loop.
We will consider several cases separately.

\smallskip

{\it Case 1: loops of small length.} The discussion
of the elementary moves above shows that loops of length 2 or
3 can be reduced to an empty loop by one elementary move.
Loops of length 1 do not exist.

\smallskip

{\it Case 2: existence of a subloop.} Assume that
the ingoing end of some $I_k$, $k\ge 1  $, coincides
with the outgoing end of some $I_l$, $k+1\le l\le n-1$.
Then the segments $I_k,\dots ,I_l$ form a subloop
of lesser length. By induction, we may assume that
it can be reduced to an empty loop by a sequence of
elementary moves. The same sequence of moves diminishes
the length of the initial loop.

\smallskip

From now on, we may and will assume that 
our loop $I_1,\dots ,I_n$ is of length $n\ge 4$
and does not contain proper subloops.

\smallskip

Applying an appropriate $g\in PSL(2, \Z )$, we may and will
assume that $\alpha =\beta =\infty$. Thus our loop starts
with $I_1=(\infty , a)$ and ends with $I_n=(b,\infty )$,
$a,b\in \Z,$ $a\ne b$ because of absence of subloops.

\smallskip

{\it Case 3: $|a-b|=1.$} We will consider the case $b=a+1$;
the other one reduces to this one by the change of orientation.
We may also assume that $I_2, \dots ,I_{n-1} \subset [a,a+1],$ 
because otherwise $\infty$ would appear once again as one of the
vertices.

\smallskip

Since $n\ge 4$, the primitive chain $I_2,\dots ,I_{n-1}$
connecting $a$ to $a+1$ has length at least 2.
We will apply to the initial loop the following
elementary moves: (i) insert $((a,a+1), (a+1,a)$ after
$I_1=(\infty ,a)$; (ii) delete $(\infty , a),((a,a+1), (a+1, \infty )$.
The resulting loop connecting $a$ to $a$
has length $n-1$ and can be reduced to the empty loop by the
inductive assumption.

\smallskip

{\it Case 4: $|a-b|\ge 2.$} Again, we may assume that
$a<b$ and that $I_2, \dots ,I_{n-1} \subset [a,b].$ 
Consider two subcases.

\smallskip

It can be that the Lebesgue measure of $I_2$  is 1, so that $I_2=(a,a+1).$
Then we apply two elementary moves: (i) insert $(a+1,\infty ),
(\infty ,a+1)$ after $I_2$; (ii) delete $(\infty ,a)=I_1, (a,a+1)=I_2,
(a+1,\infty )$. We will get a loop of length $n-1$.

\smallskip

If the Lebesgue measure of $I_2$ is $<1$, then some initial subchain
$I_2, \dots ,I_k$ of length $\ge 2$ will connect
$a$ to $a+1$. In this subcase, we will apply the following
sequence of elementary moves: (i) insert $(a,a+1), (a+1,a)$ after
$I_1=(\infty ,a)$; (ii) insert $(a+1,\infty ),
(\infty ,a+1)$ after $I_k$; (iii) delete  
$(\infty , a), (a,a+1), (a+1, \infty )$.

\smallskip
The resulting loop will have length $n+1$, however, it will also
have a subloop $I_2, \dots ,I_k, (a+1,a)$ of length $\le n-1$
and $\ge 3$. The latter can be deleted by elementary moves
in view of the inductive assumption
leaving the loop of length $\le n-2.$

\smallskip

This completes the proof of (*) and of the Theorem 1.8.

\medskip

{\bf 1.9. Pseudo--measures and generalized Dedekind symbols.}
Let $V:= \{(p,q)\in \Z^2\,|\,p\ge 1 , \roman{gcd}\,(p,q)=1\}$.

\smallskip

Slightly extending the definition given in [Fu1], [Fu2], we will
call {\it a $W$--valued generalized Dedekind symbol}
any  function $D:\, V\to W$ satisfying the functional equation
$$
D(p,q)=D(p,q+p)\,.
\eqno(1.11)
$$
The symbol $D$ can be reconstructed (at least, in the absence
of 2--torsion in $W$) from
its {\it reciprocity function} $R$ defined on the subset
$V_0:=\{(p,q)\in V\,|\, q\ge 1\}$ by the equation
$$
R(p,q):=D(p,q) - D(q,-p)\,.
\eqno(1.12)
$$
From (1.11) we get a functional equation for $R$:
$$
R(p+q,q)+R(p,p+q)=R(p,q)\,.
\eqno(1.13)
$$
Fukuhara in [Fu1] considers moreover the involution
$(p,q)\mapsto (p,-q)$. If $D$ is even 
with respect to such an involution, its reciprocity
function satisfies the additional condition
$R(1,1)=0$, which together with (1.13) suffices
for existence of $D$ with reciprocity function $R$.

\smallskip

In the following, we will work with reciprocity 
functions only. 

\medskip

{\bf 1.9.1. From pseudo--measures to reciprocity functions.}
Consider the  set $\Pi$ consisting of all primitive
segments contained in $[0,1]$.

\smallskip

If $a/p<b/q$ are ends of $I\in \Pi$ written in lowest
terms with $p,q>0$, we have $(p,q)\in V_0.$
We have thus defined a map $\Pi\to V_0$.
One easily sees that it is a bijection.

\smallskip

The involution on $\Pi$ which sends $(p,q)$ to $(q,p)$
corresponds to $I\mapsto 1-I.$

\smallskip

Let $\mu$ be a $W$--valued pseudo--measure. Define
a function $R_{\mu,0}$ which in the notation
of the previous paragraph is given by
$$
R_{\mu,0}(p,q):= \mu \left(\frac{a}{p},\frac{b}{q}\right)\,.
\eqno(1.14)
$$
Furthermore, for each $n\in \Z$ put
$$
R_{\mu,n}(p,q):= \mu \left(n+\frac{a}{p},n+\frac{b}{q}\right)\,.
\eqno(1.15)
$$

\smallskip

{\bf 1.9.2. Proposition.} {\it (a) Equations
(1.14) and (1.15) determine reciprocity functions.

\smallskip

(b) We have $R_{\mu,n}(1,1)=0$ iff $\mu (n,n+1)=0.$}

\medskip

{\bf Proof.} In order to prove (a) it suffices to establish 
the equation (1.13) for $R_{\mu,n}$. When $n=0$,
this follows from (1.2) applied to the Farey triple
$$
\alpha = \dfrac{a}{p},\ \beta = \dfrac{a+b}{p+q},\ 
\gamma = \dfrac{b}{q}\,.
$$
To get the general case, one simply shifts this triple by $n$.

\smallskip

Since $R_{\mu,0}(1,1)=\mu (0,1),$  we get (b).

\medskip

{\bf 1.10. From reciprocity functions to pseudo--measures.}
Consider now any sequence of reciprocity functions
$R_n,\, n\in \Z$, and an element $\omega\in W.$
Construct from this data a $W$--valued function 
$\widetilde{\mu}$ on the set of all primitive
segments by the following prescriptions. For positively
oriented infinite segments we put:
$$
\widetilde{\mu} (-\infty ,0):=\omega\,,
\eqno(1.16)
$$
and moreover, when $n\in\Z ,\,,n\ge 1$,
$$
\widetilde{\mu} (-\infty ,n):=\omega +
R_0(1,1)+R_1(1,1) +\dots +R_{n-1}(1,1)\,,
\eqno(1.17)
$$
$$
\widetilde{\mu} (-\infty, -n):= \omega - R_{-1}(1,1)-
\dots -R_{-n}(1,1)\,,
\eqno(1.18)
$$
$$
\widetilde{\mu} (-n, \infty ):=
R_{-n}(1,1)+R_{-n+1}(1,1) +\dots +R_{-1}(1,1)-\omega\,,
\eqno(1.19)
$$
$$
\widetilde{\mu} (n, \infty ):=
-R_{n-1}(1,1)-R_{n-2}(1,1) -\dots -R_{0}(1,1)-\omega\,.
\eqno(1.20)
$$

For positively oriented finite segments we put:
$$
\widetilde{\mu}(n+\alpha ,n+\beta ):= R_n(p,q)
\eqno(1.21)
$$
if $0\le \alpha = a/p<\beta =p/q\le 1$.

\smallskip

Finally, for negatively oriented primitive segments we prescribe
the sign change, in concordance with (1.2): 
$$
\widetilde{\mu}(\beta ,\alpha )=-\widetilde{\mu} (\alpha ,\beta ).
\eqno(1.22)
$$

{\bf 1.11. Theorem.} {\it The function $\widetilde{\mu}$
is a pre--measure. Therefore it can be uniquely extended to
a pseudo--measure $\mu$.

\smallskip

The initial family $\{R_n, \omega\}$ can be reconstructed from $\mu$
with the help of (1.15) and (1.16).}

\medskip

{\bf Proof.} In view of the Theorem 1.8, it remains only to check
the relations (1.2) for primitive loops of length 3.

\smallskip

If all ends in such a loop are finite, they form a Farey
triple in $[0,1]$ or a  Farey triple shifted by some $n\in\Z.$
In this case (1.2) follows from the functional equations
for $R_n$.

\smallskip

If one end in such a loop is $\infty$, then up to an overall
change of orientation it has the form
$$
(-\infty ,n), (n,n+1), (n+1,\infty )\,.
$$
Straightforward calculations using (1.16)--(1.21)
complete the proof.

\medskip

{\bf 1.12. Example: pseudo--measures associated with holomorphic
functions vanishing at cusps.} The spaces $\P^1(\Q )$ and $\P^1(\R )$ are embedded  into the Riemann sphere
$\P^1(\C)$ using the complex
values of the same affine coordinate $z$. The infinity
point acquires one more traditional notation $i\infty$. The upper
half plane $H =\{z\,|\,\roman{Im}\,z >0\}$ is embedded
into $\P^1(\C)$ as an open subset with the boundary
$\P^1(\R )$. The metric $ds^2:=\dfrac{|dz|^2}{(\roman{Im}\,z)^2}$
has the constant negative curvature $-1$, and
the fractional linear transformations $z\mapsto gz$,
$g\in GL^+(2,\Z )$,
act upon $H$ by holomorphic isometries. 

\smallskip

Let $\Cal{O}(H)_{cusp}$ be the space of holomorphic functions $f$
on $H$ having the following property: {\it when we approach
a cusp $\alpha \in \P^1(\Q )$ along a geodesic leading to this
cusp, $|f(z)|$ vanishes faster than $l(z_0,z)^{-N}$ for all integers $N$, 
where $z_0$ is any fixed reference point in $H$,
and  $l(z_0,z)$ is the exponential of the geodesic distance 
from $z_0$ to  $z$.}

\smallskip

If $f\in \Cal{O}(H)_{cusp}$, 
the map $\P^1(\Q )^2\to \C$:
$$
(\alpha ,\beta )\mapsto \int_{\alpha}^{\beta} f(z) dz
\eqno(1.23)
$$
is well defined and satisfies (1.1), (1.2).
(Here and henceforth we always tacitly assume
that the integration path, at least in some vicinity
of $\alpha$ and $\beta$, follows a geodesic.)

\smallskip

Considering the r.h.s. of (1.23) as a linear functional
of $f$ i.e., an element of the linear dual
space $\bold{W}:= (\Cal{O}(H)_{cusp})^*$,
we get our basic $\bold{W}$--valued pseudo--measure $\mu$.
The classical constructions with modular forms introduce some additional
structures and involve passage to a quotient of $\bold{W}$,
cf. 2.6 below.

\smallskip

The real and imaginary parts of (1.23) are pseudo--measures as well.

\smallskip

More generally, one can replace the integrand in (1.23)
by a closed real analytic (or even smooth) 1--form in $H$
satisfying the same exponential vanishing condition along cusp geodesics
as above.

\medskip

{\bf 1.13. $p$--adic analogies.} The segment $[-1,1]$ is determined by the
condition $|\alpha |_{\infty}\le 1$ so that it traditionally 
is considered as an analog of the ring
of $p$--adic integers $\Z_p$ determined by 
$|\alpha |_p\le 1$. Less traditionally, we suggest
to consider the open primitive (Farey) segments of length $\le 1/2$ in 
$[-1,1]$ as analogs of residue classes $a\,\roman{mod}\,p^m$.
Both systems of subsets share the following property:
any two subsets either do not intersect, or one of them is contained
in another one. The number $m$ then corresponds to the Farey depth
of the respective segment that is, to the length of the 
continued fraction decomposition of one end.

\smallskip

The notion of a pseudo--measure is similar to that 
of $p$--adic measure: see [Ma2], sec. 8 and 9. 
This list of analogies will be continued in
the section 3.5 below.

\bigskip

\centerline{\bf 2. Modular pseudo--measures.}

\medskip

{\bf 2.1. Modular pseudo--measures.} Let $\Gamma \subset
SL(2,\Z )$ be a subgroup of finite index. In this whole
section we assume that the group of values $W$
of our pseudo--measures is a left $\Gamma$--module. The action
is denoted $\omega\mapsto g\omega $ for $g\in \Gamma$, $\omega\in W.$  

\medskip

{\bf 2.1.1. Definition.} {\it A pseudo--measure
$\mu$  with values in $W$ is called  modular with respect to $\Gamma$
if $\mu g (\alpha ,\beta )= g[\mu (\alpha ,\beta )]$
for any $(\alpha ,\beta )$ and $g\in \Gamma$. Here
$\mu g (\alpha ,\beta )$ is defined by (1.4). }

\smallskip

We denote by $M_W(\Gamma )$ the group of all such pseudo--measures.
More generally, if $\chi$ is a character of $\Gamma$
and if multiplication by $\chi (g),\,g\in \Gamma$, makes sense in $W$
(e.~g., $W$ is a module over a ring where the values of
$\chi$ lie), we denote by
$M_W(\Gamma ,\chi )$ the group of pseudo--measures $\mu$ such
that $\mu g (\alpha ,\beta )= \chi (g)\cdot g \mu (\alpha ,\beta )$. 

\medskip

{\bf 2.1.2. Example.} Let $\Gamma= SL(2,\Z).$ Then from the
Corollary 1.7 one infers that any $\Gamma$--modular pseudo--measure $\mu$ is
uniquely determined by its single value $\mu (\infty ,0)$,
because if $I=g(\infty ,0)$ is a primitive interval, then we must have
$\mu (I) = g\mu (\infty ,0)$. In particular,
$$
\mu (\infty ,\alpha )= \sum_{k=-1}^{n-1} g_k(\alpha )\mu (\infty ,0)
\eqno(2.1)
$$ 
where the matrices $g_k(\alpha )$ are defined in (1.7).

\smallskip

For a general $\Gamma$, denote by $\{h_k\}$ a system of representatives of the coset
space $\Gamma\setminus SL(2,\Z ).$ Then the similar reasoning
shows that any $\Gamma$--modular pseudo--measure $\mu$ is
uniquely determined by its values $\mu (h_k(\infty ) ,h_k(0)).$

\medskip

{\bf 2.1.3. Modularity and generalized Dedekind symbols.}
Assume that $W$ is fixed by
the total group of shifts in $SL(2,\Z )$ fixing $\infty$:
$z\mapsto z+n,$ $n\in \Z.$ 
Assume moreover that $\mu$ is modular with respect to some subgroup
(non--necessarily of finite index) containing
the total group of shifts. Then from (1.15) one sees that
that all the reciprocity functions $R_{\mu ,n}$ coincide with $R_{\mu}:=R_{\mu ,0}$,
and from (1.17)--(1.20) it follows that $R(1,1)=0$.
Hence $\mu$ is associated with a generalized Dedekind symbol.
Conversely, starting with such a symbol, we can uniquely
reconstruct a shift--invariant pseudo--measure.

\medskip

{\bf 2.2.  A description of $M_W(SL(2,\Z ))$.} Consider now the map
$$
M_W(SL(2,\Z ))\to W :\, \mu \mapsto \mu (\infty ,0)\,.
\eqno(2.2)
$$
From 2.1.1 it follows that this is an injective group
homomorphism. Its image is constrained by several conditions.

\smallskip

Firstly, because of modularity, we have
$$
(-id) \mu (\infty , 0) = \mu (\infty ,0)\, .
\eqno(2.3)
$$
Therefore, the image of (2.2) is contained in $W_+$,
the subgroup of fixed points of $-id$. Hence we may and will consider it
as a module over $PSL(2,\Z )$. The latter group is
the free product of its two subgroups
$\bold{Z}_2$ and $\bold{Z}_3$ generated respectively by the
fractional linear transformations with matrices
$$
\sigma = \left(\matrix 0& -1\\ 1& 0\endmatrix \right),\quad
\tau = \left(\matrix 0& -1\\ 1& -1\endmatrix \right)\, .
\eqno(2.4)
$$
Now, $\sigma (\infty )=0$, $\sigma (0) =\infty $.
Hence secondly,
$$
(1+\sigma )\mu (\infty ,0) =0\,.
\eqno(2.5)
$$
And finally, 
$$
\tau (0)=1,\ \tau (1) =\infty ,\ \tau (\infty )=0\,.
$$
so that the modularity implies
$$
(1+\tau +\tau^2)\mu (\infty ,0)=0\,.
\eqno(2.6)
$$
To summarize, we get the following complex
$$
0\to M_W (SL(2,\Z )) \to W_+ \to W_+\times W_+
\eqno(2.7)
$$
where the first arrow is injective, and
the last arrow is $(1+\sigma , 1+\tau +\tau^2).$

\medskip

{\bf 2.3. Theorem.} {\it The sequence (2.7) is exact.
In other words, the map $\mu\mapsto \mu (\infty ,0)$ induces
an isomorphism}
$$
M_W (SL(2,\Z )) \cong \roman{Ker}\,(1+\sigma )\cap
 \roman{Ker}\,(1+\tau +\tau^2 ) |_{W_+}\,.
\eqno(2.8)
$$ 

\smallskip

{\bf Proof.} Choose an element $\omega$ in the r.h.s.
of (2.8). Define a $W$--valued function $\widetilde{\mu}$
on primitive segments by the formula
$$
\widetilde{\mu} (g(\infty ), g(0)) := g\omega \,.
\eqno (2.9)
$$
We will check that it is a pre--measure in the sense of 1.7.1. 
This will show that any such $\omega$ comes from
a (unique) measure $\mu$. It is then easy to check that $\mu$
is in fact $SL(2,\Z )$--modular: from the formula
(1.9) and the independence of its r.h.s. from the choice
of a primitive chain $I_k=g_k(\infty ,0)$ we get using (2.9):
$$
\mu (g(\alpha ), g(\beta ))= \sum_{k=1}^n \widetilde{\mu} (g(I_k))=
\sum_{k=1}^n \widetilde{\mu} (gg_k(\infty ,0))=
$$
$$
\sum_{k=1}^n gg_k(\omega ) =
\sum_{k=1}^n g\widetilde{\mu} (I_k) = g\mu (\alpha , \beta ).
$$ 
\smallskip

The property  (1.1) holds for $\widetilde{\mu}$ in view of (2.5).
In fact, if $(\alpha ,\beta )=g(\infty ,0)$, then
$$
(\beta ,\alpha )= g(0, \infty )= g\sigma (\infty ,0)
$$
so that
$$
\widetilde{\mu} (\alpha ,\beta )+\widetilde{\mu} (\beta ,\alpha )=
g\omega + g\sigma \omega =0.
$$
Similarly, the property (1.2) follows from (2.6).
In order to deduce this, we must check that
both types of primitive loops of length 3 considered in
the proof of the Theorem 1.8 can be represented in the form
$$
g(\infty ,0),\,g\tau (\infty ,0),\, g\tau^2 (\infty ,0)
$$
for an appropriate $g$. We leave this as an easy exercise.

\medskip

{\bf 2.3.1. Rational period functions as pseudo--measures.}
Fix an integer $k\ge 0$ and denote by $W_k$ the space of 
higher differentials $\omega (z):=q(z) (dz)^k$ where $q(z)$ is a
rational function. Define the left action of $PSL(2, \Z )$
on  $W_k$ by $(g\omega )(z) := \omega (g^{-1}(z)).$
Then the right hand side of (2.8) consists of such $q(z) (dz)^k$
that $q(z)$ satisfies the equations
$$
q(z)+z^{-2k}q(\frac{-1}{z})=0, \quad 
q(z)+z^{-2k}q(1-\frac{1}{z})+(z-1)^{-2k}q(\frac{1}{1-z})=0.
$$
These equations define the space of {\it rational period functions of weight 
$2k$.} M.~Knopp introduced them in [Kn1] and  showed that such functions
can have poles only at $0$ and real quadratic irrationalities
(the latter are ``shadows'' of closed geodesics on modular
curves.) Y.~J.~Choie and D. Zagier in [ChZ] provided
a very explicit description of them. Finally, A.~Ash
studied their generalizations for arbitrary $\Gamma$,
which can be treated as pseudo--measures using the construction
described in the next subsection.

\smallskip

{\bf 2.4. Induced pseudo--measures.} By changing the group
of values $W$, we can reduce the description of $M_W(\Gamma )$
for general $\Gamma$ to the case $\Gamma =SL(2,\Z )$.

\smallskip

Concretely, given $\Gamma$ and a $\Gamma$--module $W$,
put
$$
\widehat{W}:=\roman{Hom}_{\Gamma}(PSL(2,\Z ),W)\,.
\eqno(2.10)
$$
An element $\varphi \in \widehat{W}$ is thus a map
$g\mapsto \varphi (g)\in W$ such that $\varphi (\gamma g)=
\gamma \varphi (g)$ for all $\gamma \in \Gamma$ and
$g\in SL(2,\Z ).$ Such functions form a group with
pointwise addition, and $PSL(2, \Z )$ acts on it from the left via
$(g\varphi )(\gamma ):= \varphi (\gamma g).$

\smallskip

Any $\widehat{W}$--valued $SL(2,\Z)$--modular measure $\hat{\mu}$
induces a $W$--valued $\Gamma$--modular measure:
$$
\mu (\alpha ,\beta ) := (\hat{\mu} (\alpha , \beta )) (1_{SL(2,\Z)})\,.
\eqno(2.11)
$$
Conversely, any $W$--valued $\Gamma$--modular measure $\mu$
induces a  $\widehat{W}$--valued $SL(2,\Z)$--modular measure $\hat{\mu}$:
$$
(\tilde{\mu}(\alpha ,\beta ))(g):=\mu (g(\alpha ),g(\beta )) \,.
\eqno(2.12)
$$

\smallskip

{\bf 2.4.1. Proposition.} {\it The maps (2.11), (2.12) are well defined
and mutually inverse. Thus, they produce a canonical
isomorphism}
$$
M_W(\Gamma )\cong M_{\widehat{W}} (SL(2,\Z )).
\eqno(2.13)
$$

\smallskip

The proof is straightforward.

\medskip

{\bf 2.4.2. Modular pseudo--measures and cohomology.} Given 
$\mu \in M_W(\Gamma )$
and $\alpha \in \P^1(\Q )$, consider the function
$$
c_{\alpha}^{\mu} =c_{\alpha}:\, \Gamma \to W,\quad c_{\alpha}(g):=\mu (g\alpha , \alpha ).
$$
From (1.1), (1.2) and the modularity of $\mu$ it follows that this 
is an 1--cocycle in $Z^1(\Gamma ,W)$:
$$
c_{\alpha}(gh)=c_{\alpha}(g)+gc_{\alpha}(h)\,.
$$
Changing $\alpha$, we get a cohomologous cocycle:
$$
c_{\alpha}(g)-c_{\beta}(g) = g\mu (\alpha ,\beta )-\mu (\alpha ,\beta )\,.
$$
If we restrict $c_{\alpha}$ upon the subgroup $\Gamma_{\alpha}$
fixing $\alpha$, we get the trivial cocycle, so the respective cohomology
class vanishes. Call a cohomology class in $H^1(\Gamma ,W)$
{\it cuspidal} if it vanishes after restriction on each $\Gamma_{\alpha}.$

\smallskip

Thus, we have a canonical map of $\Gamma$--modular 
pseudo--measures to the cuspidal cohomology
$$
M_W(\Gamma )\to H^1(\Gamma ,W )_{cusp}\, .
$$

\smallskip

If $\Gamma =PSL(2,\Z )$, both groups have more compact descriptions.
In particular, the map
$$
Z^1(PSL(2,\Z ),W_+)\mapsto \roman{Ker}\, (1+\sigma )\times 
\roman{Ker}\, (1+\tau +\tau^2) \subset W_+\times W_+:\
c\mapsto (c(\sigma ),c(\tau ))
$$
is a bijection.

The value $\mu (\infty ,0) =-c_{\infty} (\sigma )$ 
furnishes the connection between the two descriptions
which takes the following form:

\medskip

{\bf 2.4.3. Proposition.} {\it (i) For any $PSL(2,\Z )$--modular 
pseudo--measure $\mu$, we have
$$
c_{\infty}^{\mu} (\sigma ) =c_{\infty}^{\mu}(\tau )=-\mu (\infty ,0)\in W_+\,.
$$ 
\smallskip

(ii) This correspondence defines a bijection
$$
M_W(PSL(2,\Z)) \cong Z^1(PSL(2,\Z ),W_+)_{cusp}
$$
where the latter group by definition consists of cocycles 
with equal components 
$c(\sigma )=c(\tau )$.}

\smallskip

In the sec. 4.8.2,
this picture will be generalized to the non--commutative case. 

\medskip

{\bf 2.5. Action of the Hecke operators on pseudo--measures.}
To define the action of the Hecke operators upon $M_W(\Gamma )$, 
we have to assume that
the group of values $W$ is a left module not only over $\Gamma$ 
(as up to now) but in fact over $GL^+(2,\Q )$. We adopt this assumption
till the end of this section.

\smallskip

Then $M_W(\Gamma )$ becomes  a  $GL^+(2,\Q )$--bimodule, with the
right action (1.4) and the left one
$$
(g\mu ) (\alpha ,\beta ):= g [\mu (\alpha ,\beta )]\, .
\eqno(2.14)
$$

\smallskip

Let $\Delta$ be a double coset in 
$\Gamma \setminus GL^+(2,\Q )/\Gamma$, or a finite union
of such cosets. Denote by $\{\delta_i\}$
a complete finite family of representatives of $\Gamma \setminus \Delta$.

\smallskip

{\bf 2.5.2. Proposition.} {\it The map
$$
T_{\Delta}:\, \mu \mapsto \mu_{\Delta}:= \sum_i \delta_i^{-1}\mu\delta_i
\eqno(2.15)
$$
restricted to $M_W(\Gamma )$ depends only on
$\Delta$ and sends  $M_W(\Gamma )$ to itself.}

\smallskip

{\bf Proof.} If $\{\delta_i\}$ are replaced by $\{g_i\delta_i\}$,
$g_i\in \Gamma$, the r.h.s. of (2.15) gets replaced by
$$
\sum_i \delta_i^{-1}g_i^{-1}\mu g_i\delta_i \,.
$$
But $g_i^{-1}\mu g_i =\mu$ on $M_W(\Gamma )$. 

\smallskip

To check that $\mu_{\Delta}g =g \mu_{\Delta}$ for $g\in \Gamma$,
notice that the right action of $\Gamma$ induces a
permutation of the set $\Gamma \setminus \Delta$ so that for
each $g\in \Gamma$,
$$
\delta_i g =g^{\prime}(i,g)\cdot \delta_{j(i,g)},\quad
\quad  g^{\prime}(i,g)\in \Gamma ,
$$
and $i\mapsto j(i,g)$ is a permutation of indices.

\smallskip

Hence
$$
\delta_i^{-1}\mu\delta_i g= 
\delta_i^{-1}\mu g^{\prime}(i,g)\cdot \delta_{j(i,g)} =
\delta_i^{-1} g^{\prime}(i,g)\mu \cdot \delta_{j(i,g)}.
$$
But $\delta_i^{-1} g^{\prime}(i,g) = g \delta_{j(i,g)}^{-1}$.
Therefore,  
$$
\delta_i^{-1}\mu\delta_i g =g\delta_{j(i,g)}^{-1}\mu \delta_{j(i,g)},
$$
and while $g$ is kept fixed, the summation over $i$ produces the same result as summation over $j(i,g)$. This completes the proof.

\smallskip

Notice in conclusion that all elements of one double coset 
$\Delta$ have the same
determinant, say, $D$. Usually one normalizes
Hecke operators by choosing $\Delta \subset M(2,\Z )$
and replacing $T_{\Delta}$ by $D T_\Delta$.

\medskip

{\bf 2.5.3. Classical Hecke operators $T_n$.} They correspond to
the case $\Gamma =SL(2, \Z )$. Let $\Delta_n$
be the finite union of the double classes
contained in $M_n(2,\Z )$: matrices with
integer entries and determinant $n$. It is well known
that the following matrices form a complete system of representatives
of $\Gamma\setminus \Delta_n$:
$$
\left(\matrix a & b\\0 & d 
\endmatrix\right) \,,\quad ad=n,\ 1\le b\le d \, .
\eqno(2.16)
$$
We put $T_n:=T_{\Delta_n}.$

\medskip

{\bf 2.6. Pseudo--measures associated with cusp forms.} 
The classical theory of modular symbols associated with
cusp forms of arbitrary weight with respect to a group
$\Gamma$ furnishes the basic examples of modular pseudo--measures.
They are specializations of the general construction of sec. 1.12.

\smallskip

We will recall the main formulas and notation.

\smallskip

First of all, one easily sees that the space 
$\Cal{O}(H)_{cusp}$ is a $\C[z]$--module,
and simultaneously a right $GL^+(2,\Z )$--module
with respect to the ``variable change'' action.
The change from the right to the left comes as a result
of dualization, cf. the last paragraph of sec. 1.

\smallskip

In the classical theory, one restricts
the pseudo--measure (1.23) onto subspaces
of $\Cal{O}(H)_{cusp}$ consisting of
the products $f(z)P(z)$ where $f$ is a classical
cusp form  of a weight $w+2$ and $P$ a polynomial
of degree $w$. The action of $GL^+(2,\Z )$
is rather arbitrarily redistributed between $f$ and $P$.

\smallskip

To present this part of the structure more
systematically, we have to consider four
classes of  linear representations (here understood as right actions) 
of  $GL^+(2,\Q )$: (i) the one--dimensional
determinantal representation $g\mapsto \roman{det}\,g\cdot id$;
(ii) the symmetric power of the basic two--dimensional representation;
(iii) the ``variable change'' action upon $\Cal{O}(H)_{cusp}$,
that is, the inverse image $g^*$ with respect to the fractional
linear action of $PGL^+(2,\Q )$
$g:\,H\to H$, $g\mapsto g(z)$; (iv) the similar inverse image map
on polydifferentials $\Omega_1(H)^{\otimes r}_{cusp}$,
that is, holomorphic tensors $f(z) (dz)^r$, $r\in \Z$;

\smallskip

The latter action is traditionally translated into a
``higher weight'' action on functions $f(z)$
from $\Cal{O}(H)_{cusp}$ via dividing by $(dz)^r$,
and further tensoring by a power of the determinantal
representation. As a result, the picture becomes
somewhat messy, and moreover, in different expositions
different normalizations and distributions of determinants
between $f$ and $P$ are adopted. Anyway, we will fix
our choices by the following conventions (the same as in [He1], [He2].)

\smallskip

For an integer $w\ge 0$, define the right action
of weight $w+2$ upon holomorphic (or meromorphic)
functions on $H$ by
$$
f|[g]_{w+2}(z):= (\roman{det}\, g)^{w+1}
f(gz) j(g,z)^{-(w+2)}
\eqno(2.17)
$$
where we routinely denote $j(g,z):=cz+d$ for
$$
g=\left(\matrix a & b\\c & d 
\endmatrix\right) \, .
$$

Moreover, define the right action of $GL(2, \R)$
on polynomials in two variables by
$$
(Pg)(X,Y):=  P( (\roman{det}\, g)^{-1}(aX+bY),(\roman{det}\, g)^{-1}(cX+dY))\, .
\eqno(2.18)
$$
Let now $\alpha$, $\beta$ be two points in $H\cup \P^1(\Q )$.
Then we have, for any homogeneous
polynomial $P(X,Y)$ of degree $w$ and $g\in GL^+(2,\Q )$, 
in view of (2.17) and (2.18),
$$
\int_{g\alpha}^{g\beta} f(z)P(z,1)dz=
\int_{\alpha}^{\beta} f(gz)P(gz,1)d(gz) =
$$
$$
\int_{\alpha}^{\beta} f|[g]_{w+2}(z)
\cdot  (\roman{det}\, g)^{-w-1}\cdot j(g,z)^{w+2}\cdot
P\left(\frac{az+b}{cz+d},1\right)\cdot j(g,z)^{-2} \cdot\roman{det}\, g\cdot dz =
$$
$$
\int_{\alpha}^{\beta} f|[g]_{w+2}(z)
P( (\roman{det}\, g)^{-1}(az+b),(\roman{det}\, g)^{-1}(cz+d)) dz=
$$
$$
\int_{\alpha}^{\beta} f|[g]_{w+2}(z) (Pg)(z,1) dz.
\eqno(2.19)
$$
In particular, if for given $f$, $g$, and a constant $\varepsilon$ 
we have
$$
f|[g]_{w+2}(z) = \varepsilon f(z),
\eqno(2.20)
$$ 
then
$$
\int_{g\alpha}^{g\beta} f(z)P(z,1)dz = 
\int_{\alpha}^{\beta} f(z) \varepsilon (Pg)(z,1) dz.
\eqno(2.21)
$$
More generally, if for a finite family of $g_k\in GL^+(2,\Q )$,
$c_k\in \C$,
we have
$$
\sum_k c_kf|[g_k]_{w+2}(z) = \varepsilon f(z),
\eqno(2.22)
$$ 
then
$$
\sum_k c_k\int_{g_k\alpha}^{g_k\beta} f(z)P(z,1)dz = 
\int_{\alpha}^{\beta} f(z) \varepsilon \sum_k c_k (Pg_k)(z,1) dz.
\eqno(2.23)
$$
These equations are especially useful in the standard context of
modular and cusp forms and Hecke operators. 

\smallskip

Let $\Gamma \subset
SL(2,\Z )$ be a subgroup of finite index, $w\ge 0$ an integer.
Recall that a cusp form $f(z)$ of weight $w+2$ with respect to $\Gamma$
is a holomorphic function on the upper half--plane $H$,
vanishing at cusps and such that for all $g\in \Gamma$,
(2.20) holds with constant $\varepsilon =1$. More generally,
we can consider constants $\varepsilon =\chi (g)$,
where $\chi$ is a character of $\Gamma$.

\smallskip

When $w$ is odd, nonvanishing cusp form can exist only
if $-I\notin \Gamma$ where $I$ is the identity matrix.
In particular, for $\Gamma =SL(2,\Z )$ there can exist
only cusp forms of even weight.

\smallskip

Denote by $S_{w+2}(\Gamma )$ (resp. $S_{w+2}(\Gamma ,\chi ))$ the complex space of cusp forms of weight $w+2$ for $\Gamma$ (resp. cusp forms 
with character $\chi$). Let $F_w$ be the space of polynomial forms
of degree $w$ in two variables. Let $W$ be the space
of linear functionals on $S_{w+2}(\Gamma )\otimes F_w.$
Then the function on $(\alpha ,\beta )\in {\P}^1(\Q )^2$
with values in the space $W$ 
$$
\mu (\alpha ,\beta ):\, f\otimes P \mapsto \int_{\alpha}^{\beta} f(z) P(z,1) dz
\eqno(2.24)
$$
is a pseudo--measure. 
\smallskip

We can call this pseudo--measure ``the shadow'' of the respective
modular symbol. 

\smallskip

Formulas (2.17) and (2.18) define the structure of a right
$\Gamma$--module upon $S_{w+2}(\Gamma )\otimes F_w$
and hence the dual structure of a left $\Gamma$--module upon $W$.   
Formula (2.21) then shows that 
that the pseudo--measure (2.24) is modular with respect to $\Gamma$.

\bigskip

\centerline{\bf 3. The L\'evy functions and the L\'evy--Mellin transform}

\medskip

{\bf 3.1. The L\'evy functions.} A classical L\'evy function
(see [L]) $L(f)(\alpha )$ of a real argument $\alpha$ is given by
the formula 
$$
L(f)(\alpha ):=\sum_{n=0}^{\infty} f(q_n(\alpha ),q_{n+1}(\alpha ))
\eqno(3.1)
$$
where $q_n(\alpha ),$ $n\ge 0$, is the sequence of denominators 
of normalized convergents to $\alpha$ (see (1.5)), and $f$ is a function
defined on pairs $(q^{\prime},q) \in \Z^2$, $1\le q^{\prime}\le q$,
$\roman{gcd}\,(q^{\prime},q)=1$, taking values in a topological
group and sufficiently quickly decreasing so that (3.1)
absolutely converges. Then $L(f)(\alpha )$ is continuous
on irrational numbers. Moreover, it has period 1.

\smallskip

It was remarked in [MaMar1] that for certain simple $f$
related to modular symbols, the integral $\int_0^1 L(f)(\alpha )d\alpha$
is a Dirichlet series directly related to
the Mellin transform of an appropriate cusp form. 

\smallskip

In this section, we will develop this remark and considerably generalize
it in the context of modular pseudo--measures.
We will call the involved integral representations
{\it the  L\'evy--Mellin transform.} They are ``shadows'' of the 
classical Mellin transform.

\smallskip

As P.~L\'evy remarked in [L], for any $(q^{\prime},q)\ne (1,1)$
as above all $\alpha\in [0,1/2]$ such that $(q^{\prime},q)=
(q_n(\alpha ),q_{n+1}(\alpha ))$ fill a primitive semi--interval $I$
of length $(q(q+q^{\prime}))^{-1}$. In $[0,1]$, such $\alpha$
fill in addition the semi--interval $1-I$
(in [MaMar1], sec. 2.1, we have inadvertently overlooked
this symmetry).

\smallskip

Therefore the class of functions $L(f)$ (restricted to irrational numbers)
is contained in a more general
class of (formal) infinite linear combinations of 
characteristic functions of primitive segments $I$:
$$
L(f) := \sum_I f(I)\chi_I\, ,
\eqno(3.2)
$$
where $\chi_I(\alpha )=1$ for $\alpha \in I$, $0$ for  $\alpha \notin I$.

\smallskip

A family of coefficients, i.e. a map $I\mapsto f(I)$ from the set of all 
primitive segments (positively oriented, or without 
a fixed orientation) in [0,1]
to an abelian group $M$, will also be referred to as 
a L\'evy function. 

\smallskip

To treat this generalization  systematically, we will display
in the next subsection the relevant combinatorics of
primitive segments.

\medskip

{\bf 3.2. Various enumerations of primitive segments.}  
A matrix in $GL(2,\Z )$ is called {\it reduced}
if its entries are non--negative, and non--decreasing to the right
in the rows and downwards in the columns.
For a more thorough discussion, see [LewZa].

\smallskip

Clearly, each
row and column of a matrix in $GL(2,\Z )$, in particular,
reduced ones, consists of co--prime entries.
There is exactly one reduced matrix with lower row $(1,1)$.
This is one marginal case which does not quite fit
in the pattern of the following series of bijections
between several sets described below.

\smallskip

{\it The set $L$.} It consists of pairs $(c,d)$, $1\le c < d \in \Z$, $gcd\,(c,d)=1$.

\smallskip

{\it The set $R$.} It consists of pairs of reduced matrices
$(g^-,g^+)$ with one and the same lower row and determinants,
respectively, $-1$ and $+1$.

\smallskip

It is easy to see that for any pair $(c,d)\in L$,
there exists exactly one pair $(g^-_{c,d},g^+_{c,d})\in R$
with lower row $(c,d)$ so that we have a natural bijection
$L\to R$.

\smallskip

In the marginal case, only one reduced matrix $g^-_{1,1}$ exists.

\smallskip

{\it The set $S$.} It consists of pairs of primitive
segments $(I^-,I^+)$ of length $<1/2$ such that $I^-\subset [0,\dfrac{1}{2}]$,
$I^+=1-I^- \subset [\dfrac{1}{2},1]$. There is a well defined map
$R\to S$ which produces from each pair $(g^-_{c,d},g^+_{c,d})\in R$
the pair 
$$
I^-_{c,d}:= [g^-_{c,d}(0),g^-_{c,d}(1)],\quad
I^+_{c,d}:= [g^+_{c,d}(0),g^+_{c,d}(1)].
$$
Again, it is a well defined bijection.

\smallskip

In the marginal case, we get one primitive segment
$I^-_{1,1} =[0,\dfrac{1}{2}]$; it is natural to
complete it by $I^+_{1,1} :=[\dfrac{1}{2},1]$.

\smallskip

{\it The set $C$.} One element of this set
is defined as a pair $(q^{\prime},q)\in L$ such that there exists
an $\alpha$, $0<\alpha \le \dfrac{1}{2}$, for which
$(q^{\prime},q)$ is a pair of denominators $(q_n(\alpha),
q_{n+1}(\alpha ))$ of two consecutive convergents to $\alpha$.

\smallskip

The sequence of consecutive convergents to a rational
number $\alpha$ stabilizes at some $p_{n+1}/q_{n+1}$,
so only a finite number of pairs $(q^{\prime},q)$
is associated with $\alpha$.
Integers, in particular $0$ and $1$, correspond to
the marginal pair (1,1).

\smallskip

According to a lemma, used by P.~L\'evy in [L],
all such $\alpha$ fill precisely the
semi--interval $I^-_{q^{\prime},q}$, with the end
$g^-_{q^{\prime},q}(1)$ excluded.

\smallskip

Moreover, those $\alpha$ which belong to the other
half of $[0,1]$ and as well admit
$(q^{\prime},q)$ as a pair of denominators 
of two consecutive convergents, fill the
semi--interval $[g^+_{q^{\prime},q}(0),g^+_{q^{\prime},q}(1))$.
This produces one more bijection $C=L\to S$, or rather an
interpretation of a formerly constructed one.

\smallskip

The pairs of convergents involved are also encoded in this picture:
they are  $(g^-_{q^{\prime},q}(0),g^-_{q^{\prime},q}(\infty ))$
and $(g^+_{q^{\prime},q}(0),g^+_{q^{\prime},q}(\infty ))$  
respectively.

\smallskip

Returning now to L\'evy functions, we can rewrite (3.2)
as
$$
L(f)(\alpha ):= f(I_{1,1}) +\sum_{I^-_{c,d}\ni\alpha} f(I_{c,d}^-)   
+\sum_{I^+_{c,d}\ni\alpha} f(I_{c,d}^+) .
\eqno(3.3)
$$
or else as a classical L\'evy function
$$
L(f)(\alpha ):= f(1,1) +\sum_{n\ge 0} f^{-}(q_n(\alpha ),q_{n+1}(\alpha ))
\chi_{(0,1/2]} (\alpha ) +
$$
$$ 
\sum_{n\ge 0} f^{+}(q_n(\alpha ),q_{n+1}(\alpha ))
\chi_{[1/2,1)} (\alpha )  
\eqno(3.4)
$$
if the usual proviso about convergence is satisfied.

\smallskip

Notice that for any $(c,d)\in L$, the number $n$ such that
$(c,d)=(q_n(\alpha ),q_{n+1}(\alpha ))$ is uniquely 
defined by $(c,d)$ and a choice of sign
$\pm$ determining the position of $\alpha$ in the right/left
half of $[0,1]$. The sequence of convergents preceding the $(n+1)$--th
one is determined as well. This means that choosing 
$f^{\pm}$, we may explicitly refer to all these data,
including the sequence of incomplete quotients up  to the $(n+1)$--th
one. Cf. examples in sec. 2.1 and 2.2.1 of [MaMar1].

\medskip

{\bf 3.3. The L\'evy--Mellin transform.} Below we will use formal Dirichlet series.

\smallskip

Let $\Cal{A}$ be an abelian group. With any sequence
$A=\{a_1,\dots ,a_n, \dots\}$  of elements of $\Cal{A}$ 
we associate a formal expression 
$$
L_A(s):=\sum_{n=1}^{\infty} \frac{a_n}{n^s} \,.
$$
If $\Cal{B},\Cal{C}$ are  abelian groups,
and we are given a composition $\Cal{A}\times \Cal{B}\to\Cal{C}:$
$(a,b)\mapsto a\cdot b$ which is a group homomorphism,
we have the induced composition
$$
(a_1,\dots ,a_n, \dots )\cdot (b_1,\dots ,b_n, \dots )=
(c_1,\dots ,c_n, \dots ), \quad c_n:=\sum_{d_1,d_2:\,d_1d_2=n}a_{d_1}\cdot b_{d_2}
$$
which is interpreted as the  multiplication of the respective
formal Dirichlet series:
$$
L_C(s) = L_A(s)\cdot L_B(s).
$$
In particular, we can multiply any $L_A(s)$ by a Dirichlet series
with integer coefficients.

\smallskip

Finally, the operator $\{a_1,\dots ,a_n, \dots\}\mapsto
\{1^wa_1, 2^wa_2 ,\dots ,n^wa_n, \dots\}$ corresponds to
the argument shift $L_A(s-w)$. Here generally $w\in \Z$, $w\ge 0$,
but one can use also negative $w$ if $\Cal{A}$ is a linear $\Q$--space.

\smallskip

Now let $W$ be a left $GL^+(2,\Q )$--module: $w\mapsto g[w]$.
We will apply this formalism to the Dirichlet series
with coefficients in $\Z [GL^+(2,\Q )]$, in $W$, and the
composition induced by the above action.

\smallskip

Let $\mu \in M_W(SL(2, \Z ))$
(for notation, see sec. 2.1.1.) Consider the following L\'evy function
$f_{\mu}$ non--vanishing only on the primitive segments in $[0,1/2]$:
$$
f_{\mu}(I^-_{c,d})(s):= \frac{1}{|I^-_{c,d}| d^{s}}\, 
\left(\matrix 1 & -cd^{-1}\\0 & d^{-1} 
\endmatrix\right) \, 
\left[\mu(\infty , \frac{c}{d})\right] \,,
\eqno(3.5)
$$
where the matrix preceding the square brackets acts upon 
$\mu(\infty , \frac{c}{d})$.
This L\'evy function takes values in the group of formal Dirichlet series
with coefficients in $W$, with only one non--vanishing 
term of the series.

\smallskip

{\bf 3.3.1. Definition.} {\it The Mellin--L\'evy transform
of the pseudo--measure $\mu$ is the formal Dirichlet series with
coefficients in $W$:}
$$
LM_{\mu}(s):=\int_0^{1/2} L(f_{\mu})(\alpha ,s) d\alpha \,.
\eqno(3.6)
$$
(The formal argument $s$ is included in the notation
for future use.)

\smallskip

Moreover, introduce the following two formal Dirichlet series
with coefficients in $\Z [GL^+(2,\Z )]$:
$$
Z_-(s):=\sum_{d_1=1}^{\infty} \left(\matrix 1 & 0\\0 & d_1^{-1} 
\endmatrix\right) \,\frac{1}{d_1^s}\,,\quad
Z_+(s):=\sum_{d_2=1}^{\infty} \left(\matrix d_2^{-1} & 0\\0 & 1  
\endmatrix\right) \,\frac{1}{d_2^s}\, .
$$

\smallskip

{\bf 3.4. Theorem.} {\it We have the following identity
between the formal Dirichlet series:
$$
Z_+(s)\cdot Z_-(s)\cdot LM_{\mu}(s)=\sum_{n=1}^{\infty}\frac{(T_n\mu )(\infty ,0)}{n^s}
\eqno(3.7)
$$
where $T_n$ is the Hecke operator described in 2.5.3 and 2.5.2.} 

\smallskip

{\bf Proof.} Put  $\overline{L}: = L\cup \{(1,1)\}$,
where the set $L$ was defined in sec. 3.2. Each matrix from (2.16)
representing a coset in $SL(2,\Z )\setminus M_n(2,\Z )$
can be uniquely written as
$$
\left(\matrix d_2 & cd_1\\0 & dd_1 
\endmatrix\right)
\eqno(3.8)
$$
where $d, d_1, d_2$ are natural numbers with $dd_1d_2=n$
and $(c,d)\in \overline{L}.$ Acting from the left
upon $(\infty ,0)$ this matrix produces $(\infty ,cd^{-1}).$
In the respective summand of $(T_n\mu )(\infty ,0)$
(cf. (2.15)), the inverse matrix to (3.8)
acts on the left. We have
$$
\left(\matrix d_2 & cd_1\\0 & dd_1 
\endmatrix\right)^{-1}=
\left(\matrix d_2^{-1} & 0\\0 &1 
\endmatrix\right)\cdot
\left(\matrix 1 & 0\\0 & d_1^{-1} 
\endmatrix\right)\cdot
\left(\matrix 1 & -cd^{-1}\\0 & d^{-1} 
\endmatrix\right)\,\,.
\eqno(3.9)
$$
Summing over all $d,d_1,d_2, c$ we get the identity (3.7)
between the formal Dirichlet series.

\smallskip

This Theorem  justifies the name of the Mellin--L\'evy 
transform. In fact, applying (3.7) to a pseudo--measure
associated with the space of $SL(2,\Z )$--cusp forms
of a fixed weight, and taking the value
of the r.h.s. of (3.7) on an eigen--form
for all Hecke operators, we get essentially the usual Mellin
transform of this form. The l.h.s then furnishes
its representation as an integral over a real segment
(up to extra $Z_{\pm}$--factors) replacing the
more common integral along the upper imaginary
half--line.

\smallskip

Of course, one can establish versions of this theorem for
the standard congruence subgroups $\Gamma_1(N)$,
but (3.7) is universal in the same sense as the Fourier
expansion presented in [Me] is universal.

\medskip

{\bf 3.5. $p$--adic analogies continued.} Returning to the
discussion in sec. 1.13, we now suggest the reader to compare
the formula (3.5) with the definition of the $p$--adic measure
given by the formula (39) of [Ma2].

\smallskip
 
We hope that this list makes convincing our suggestion
that the theory of pseudo--measures can be considered as
an $\infty$--adic phenomenon.

\bigskip

\centerline{\bf 4. Non--commutative pseudo--measures}

\centerline{\bf and iterated L\'evy--Mellin transform}

\medskip

 {\bf 4.1. Definition.}  {\it A pseudo--measure on
${\P}^1(\R )$ with values in a non (necessarily) commutative group  (written 
multiplicatively) 
$U$ is a function $J :\, {\P}^1(\Q )\times {\P}^1(\Q )\to U$,
$(\alpha ,\beta ) \mapsto J_{\alpha}^{\beta} \in U$
satisfying the following conditions: for any $\alpha ,\beta ,\gamma \in
{\P}^1(\Q )$},
$$
J_{\alpha}^{\alpha} =1,\quad J_{\beta}^{\alpha}J_{\alpha}^{\beta}=1,\quad
J^{\alpha}_{\gamma}J^{\gamma}_{\beta}J_{\alpha}^{\beta}=1\, .
\eqno(4.1)
$$

\smallskip

The formalism of the sections 1 and 2 can be partially generalized
to the non--commutative case.

\medskip

{\bf 4.1.1. Identical and inverse pseudo--measures.} The set $M_U$ of
$U$--valued pseudo--measures generally does not form a group,
but only a set with involution and
a marked point: $(J^{-1})_{\alpha}^{\beta} := (J^{\alpha}_{\beta})^{-1}$
and $J_{\alpha}^{\beta} \equiv 1_U$  are pseudo--measures.

\medskip

{\bf 4.1.2. Universal pseudo--measure.} Let $\bold{U}$
be a free group freely generated by the set $\Q$. 
Let $\langle \alpha\rangle \in \bold{U}$ be the generator
corresponding to $\alpha$. The map
$$
(\alpha ,\beta )\mapsto \langle \beta\rangle\langle \alpha\rangle^{-1},
\quad (\infty ,\beta )\mapsto  \langle \beta\rangle\,,
$$
is a pseudo--measure. It is universal in an evident sense
(cf. sec. 1.4.)

\medskip

{\bf 4.1.3. Primitive chains.} Let $I_j=(\alpha_j ,\alpha_{j+1})$,
$j=1, \dots ,n$,
be a primitive chain connecting $\alpha :=\alpha_1$
to $\beta := \alpha_{n+1}$ as at the end of sec. 1.6. Then
$$
J_{\alpha}^{\beta} = J_{\alpha_n}^{\alpha_{n+1}} J_{\alpha_{n-1}}^{\alpha_n}
\dots J_{\alpha_1}^{\alpha_2} \, .
\eqno(4.2)
$$
In particular, each pseudo--measure is determined by its values
on primitive segments.

\medskip

{\bf 4.2. Definition.} {\it An $U$--valued pre--measure $\widetilde{J}$
is a function on primitive segments satisfying the relations
(4.1) written for primitive chains only.}

\medskip

{\bf 4.3. Theorem.} {Each pre--measure $\widetilde{J}$ can be uniquely extended
to a pseudo--measure $J$.}

\medskip

The proof proceeds exactly as that of the Theorem 1.8, with only minor
local modifications. We define $J$ by any of the formulas
$$
J_{\alpha}^{\beta} = \widetilde{J}_{\alpha_n}^{\alpha_{n+1}} \widetilde{J}_{\alpha_{n-1}}^{\alpha_n}
\dots \widetilde{J}_{\alpha_1}^{\alpha_2} \, .
\eqno(4.3)
$$
using a primitive chain as in 4.1.3, and then prove that this
prescription does not depend on the choice of this
chain using the argument with elementary moves.
One should check only that elementary moves are compatible
with non--commutative relations (4.1) which is evident.

\medskip

{\bf 4.4. Non--commutative reciprocity functions.} 
By analogy with 1.9, we can introduce the notion
of a non--commutative $U$--valued reciprocity function
$R_p^q$. The extension of the functional equation (1.13)
should read
$$
R_{p+q}^qR_p^{p+q}=R_p^q\,.
\eqno(4.4)
$$
The analysis and results of 1.9--1.11 (not involving Dirichlet symbols)
can be easily transported to this context.

\medskip

{\bf 4.5. Right action of $GL^+(2,\Q ).$} This group acts
upon the set of pseudo--measures $M_U$ from the right as in the commutative case:
$$
(Jg)_{\alpha}^{\beta} := J_{g\alpha}^{g\beta} \,.
$$

\medskip

{\bf 4.6. Modular pseudo--measures.} Assume now that a subgroup
of finite index $\Gamma \subset SL(2, \Z )$ acts upon $U$
from the left by group automorphisms, $u\mapsto gu.$ 
As in the commutative case,
we call an $U$--valued pseudo--measure $J$ {\it modular with respect
to $\Gamma$}
if it satisfies the condition
$$
(Jg)_{\alpha}^{\beta} = g[J_{\alpha}^{\beta}]
$$
for all $g\in \Gamma$, $\alpha ,\beta \in \P^1(\Q )$.
Denote by $M_U(\Gamma )$ the pointed set of such measures.
We can repeat the argument and the construction of
2.1.2 showing that any element $J\in M_U(\Gamma )$
is uniquely determined by the values $J_{h_k(\infty )}^{h_k(0)}\in U$ 
where $\{h_k\}$ runs over a system of representatives
of $\Gamma \setminus SL(2,\Z ).$

\smallskip

An analog of 2.2 and of  the Theorem 2.3 holds as well.

\smallskip

Namely, consider the map
$$
M_U(SL(2,\Z ))\to U :\, J \mapsto J_{\infty}^0\,.
$$
From the argument above it follows that this is an injective map
of pointed sets with involution. Its image is constrained by the 
similar conditions as in the abelian case.

\smallskip

Firstly, because of modularity, we have
$$
(-id) [J_{\infty}^0] = J_{\infty}^0\, .
$$
Therefore, the image of (3.1) is contained in $U_+$,
$(-id)$--invariant subgroup of $U$. Hence we may and will consider 
this subgroup
as a (non--commutative) module over $PSL(2,\Z )$. 
Secondly,
$$
\sigma[J_{\infty}^0]\cdot J_{\infty}^0 =1_U\,.
\eqno(4.5)
$$
And finally, 
$$
\tau^2[J_{\infty}^0]\cdot \tau[J_{\infty}^0]\cdot J_{\infty}^0=1_U\,.
\eqno(4.6)
$$
To summarize, we get the following morphisms of pointed sets
$$
M_U (SL(2,\Z )) \to U_+ \to U_+\times U_+
\eqno(4.7)
$$
where the first  arrow is the embedding described above
and the second one, say $\varphi$, sends $u$ to $(\sigma u\cdot u ,$ $
 \tau^2u\cdot\tau u\cdot u).$

\medskip

{\bf 4.7. Theorem.} {\it The sequence of pointed sets (4.7) is exact.
In other words, the map $J\mapsto J_{\infty}^0$ induces
a bijection}
$$
M_U (SL(2,\Z )) \cong \varphi^{-1} (1_U,1_U)
\eqno(4.8)
$$ 
\smallskip

{\bf Sketch of proof.} The proof follows exactly the same plan as that
of Theorem 2.3. The only difference is that now
the relations (1.1), (1.2) are replaced by their
non--commutative versions (4.1), and the prescription
(2.1) must be replaced by its non--commutative version 
based upon (4.2).

\medskip

{\bf 4.8. Induced pseudo--measures.} As in the commutative case,
by changing the group
of values $U$, we can reduce the description of $M_U(\Gamma )$
for general $\Gamma$ to the case $\Gamma =SL(2,\Z )$.

\smallskip

Namely, put 
$$
\widehat{U}:=\roman{Map}_{\Gamma}(PSL(2,\Z ),U)\,.
\eqno(4.9)
$$
An element $\varphi \in \widehat{U}$ is thus a map
$g\mapsto \varphi (g)\in U$ such that $\varphi (\gamma g)=
\gamma \varphi (g)$ for all $\gamma \in \Gamma$ and
$g\in SL(2,\Z ).$ Such functions form a group with pointwise
multiplication, 
and $PSL(2, \Z )$ acts on it via
$(g\varphi )(\gamma ):= \varphi (\gamma g).$

\smallskip

Any $\widehat{U}$--valued $SL(2,\Z)$--modular measure $\widehat{J}$
induces an $U$--valued $\Gamma$--modular measure $J$:
$$
J_{\alpha}^{\beta}  := \widehat{J}_{\alpha}^{\beta} (1_{SL(2,\Z)})\,.
\eqno(4.10)
$$
Conversely, any $U$--valued $\Gamma$--modular measure $J$
produces a  $\widehat{U}$--valued $SL(2,\Z)$--modular measure $\widehat{J}$:
$$
\widehat{J}_{\alpha}^{\beta} (g):=J_{g(\alpha )}^{g(\beta )} \,.
\eqno(4.11)
$$

{\bf 4.8.1. Proposition.} {\it The maps (4.10), (4.11) are well defined
and mutually inverse. Thus, they produce a canonical
isomorphism}
$$
M_U(\Gamma )\cong M_{\widehat{U}} (SL(2,\Z )).
\eqno(4.12)
$$

\smallskip

The proof is straightforward.

\medskip

{\bf 4.8.2. Pseudo--measures and cohomology.} As in the
commutative case, given $J \in M_U(\Gamma )$
and $\alpha \in \P^1(\Q )$, consider the function
$$
c_{\alpha}^J=c_{\alpha}:\, \Gamma \to U,\quad c_{\alpha}(g):=J^{\alpha}_{g\alpha}.
$$
From (4.1) and the modularity of $J$ it follows that this 
is a non--commutative 1--cocycle (we adopt the normalization
as in [Ma3]):
$$
c_{\alpha}(gh)=c_{\alpha}(g)\cdot gc_{\alpha}(h)\,.
$$
Changing $\alpha$, we get a cohomologous cocycle:
$$
c_{\beta}(g) = J_{\alpha}^{\beta} c_{\alpha}(g) 
\left(gJ_{\alpha}^{\beta}\right)^{-1}.
$$
As in the commutative case, if we restrict $c_{\alpha}$ upon the subgroup $\Gamma_{\alpha}$
fixing $\alpha$, we get the trivial cocycle, so the respective cohomology
class vanishes. This furnishes a canonical map of modular measures to the cuspidal cohomology
$$
M_U(\Gamma )\to H^1(\Gamma ,U )_{cusp}\, .
\eqno(4.13)
$$

\smallskip

Again, for $\Gamma =PSL(2,\Z )$ we have two 
independent descriptions of these sets,
connected by the value $\mu (\infty ,0) =c_{\infty} (\sigma )$.

\medskip

{\bf 4.8.3. Proposition.} {\it (i) For any $PSL(2,\Z )$--modular 
pseudo--measure $J$, we have
$$
c_{\infty}^{J} (\sigma ) =c_{\infty}^{J}(\tau )=J^\infty_0\in U_+\,.
$$ 
\smallskip

(ii) This correspondence defines a bijection
$$
M_W(PSL(2,\Z)) \cong Z^1(PSL(2,\Z ),U_+)_{cusp}
$$
where the latter group by definition consists of cocycles 
with equal components 
$c(\sigma )=c(\tau )$.}

\smallskip

This easily follows along the lines of [Ma4], Proposition 1.2.1.

\medskip

{\bf 4.9. Iterated integrals I.} In this and the next subsection we
describe a systematic way to construct  non--commutative pseudo--measures
$J$  based upon iterated integration.
We start with the classical case, when $J$ is obtained
by integration along geodesics connecting cusps. This construction
generalizes  that of sec. 1.12 above, which furnishes
the ``linear approximation'' to it.

\smallskip

Consider a pseudo--measure $\mu :\P^1(\Q)^2\to W$ 
given by the formulas (1.23)
restricted to a finite--dimensional subspace $W^*\subset \Cal{O}(H)_{cusp}$,
that is, induced by the respective finite--dimensional
quotient of $\bold{W}$ from the last paragraph of sec. 1.

\smallskip

Consider the ring of formal non--commutative series 
$\C \langle\langle W\rangle\rangle$ which is the completion
of the tensor algebra of $W$ modulo powers of the augmentation ideal $(W)$. 
Put $U=1+(W)$. This is a multiplicative subgroup
in this ring which will be the group of values of the following
pseudo--measure $J$:
$$
J_{\alpha}^{\beta}:= 1+\sum_{n=1}^{\infty} J_{\alpha}^{\beta}(n)\,.
\eqno(4.14)
$$
Here $J_{\alpha}^{\beta}(n)\in W^{\otimes n}$ is defined
as the linear functional upon $W^{*\otimes n}$ whose value at
$f_1\otimes \dots \otimes f_n$, $f_i\in W^*$, is given by
the iterated integral
$$
J_{\alpha}^{\beta}(n)(f_1,\dots ,f_n):=\int_{\alpha}^{\beta} f_1(z_1)\,dz_1\int_{\alpha}^{z_1} f_2(z_2)\,dz_2\dots    \int_{\alpha}^{z_{n-1}}f(z_n)\,dz_n \,.
\eqno(4.15)
$$
The integration here is done over the  simplex $\sigma_n(\alpha ,\beta )$
consisting of the points $\beta >z_1 > z_2 > \dots z_n >\alpha$,
the sign $<$ referring to the ordering along the geodesic
oriented from $\alpha$ to $\beta$.

\smallskip

The basic properties of (4.14), including the pseudo--measure
identities (4.1), are well known, cf. a review in sec. 1 of [Ma3].
In particular, all $J_{\alpha}^{\beta}$ belong to the subgroup of
the so called {\it group--like} elements of $U$. This property
compactly encodes all the {\it shuffle relations}
between the iterated integrals (4.15).

\smallskip

Moreover,
(4.13) is functorial with respect to the variable--change
action of $GL^+(2,\Q )$ so that if the linear term
of (4.13), that is, the pseudo--measure (1.23), is $\Gamma$--modular,
then $J$ is $\Gamma$--modular as well.

\medskip

{\bf 4.10. Iterated integrals II.} We can now imitate this construction
using iterated integrals of, say, piecewise
continuous functions along segments in $\P^1(\R )$
in place of geodesics. The formalism remains exactly the same, and the general results as well; of course, in this generality it 
has nothing to do with specifics of the situation we have been
considering so far.

\smallskip

To re--introduce these specifics, we will iterate integrals 
of L\'evy functions (3.1)--(3.4) and particularly integrands
at the r.h.s. of (3.6). We will get in this way
the ``shadow'' analogs of the more classical multiple Dirichlet series
considered  in [Ma3], and shuffle relations between them.
It would be interesting to see whether one can in this way 
get new relations between the classical series. Here we only say
a few words about the structure of the resulting series.

\smallskip

Consider an iterated integral of the form (4.15) in which
$f_k(z)$ are now characteristic functions of finite segments $I_k$
in $\P^1(\R ).$ The following lemma is easy.

\medskip

{\bf 4.10.1. Lemma.} {\it The integral (4.15) as a function 
of $\alpha , \beta$ and $2n$ ends of $I_1,\dots ,I_n$ is
a piecewise polynomial form of degree $n$. 
This form depends only on the relative order of all
its arguments in $\R$.}

\medskip

Consider now a $\Gamma$--modular $W$--valued
pseudo--measure $\mu$ and a family
of L\'evy functions with coefficients in formal Dirichlet
series as in (3.5). Then we can interpret
an iterated integral of the form (4.15)
with integrands $f_{\mu}(I^-_{c_k,d_k})(s_k)\cdot \chi_{I^{-}_{c_k,d_k}}(z)$, $k=1,\dots ,n$,
as taking values directly in $W^{\otimes n}$.
Accordingly, the iterated version of the Mellin--L\'evy transform
(3.6) will be represented by a formal Dirichlet--like
series involving coefficients which are polynomial forms 
of the ends of primitive segments involved.

\bigskip

\centerline{\bf 5. Pseudo--measures and limiting modular symbols}

\medskip

{\bf 5.1. Pseudo--measures and the disconnection space.}
We recall the following definition of {\it analytic} pseudo--measures 
on totally disconnected spaces and their relation to currents 
on trees, cf. [vdP].

\medskip

{\bf 5.1.1. Definition.} {\it 
 Let $\Omega$ be a totally disconnected compact Hausdorff space and
let $W$ be an abelian group. Let $C(\Omega,W)=C(\Omega,\Z)\otimes_\Z W$ 
denote the group of locally constant $W$--valued continuous functions 
on $\Omega$. An analytic $W$--valued pseudo--measure on $\Omega$ is a 
map $\mu: C(\Omega,W)\to W$ satisfying the properties:

\smallskip

(i) $\mu(V \cup V')=\mu(V)+\mu(V')$ for $V,V'\subset \Omega$ clopen 
subsets with $V\cap V'=\emptyset$, where we identify a set with its
characteristic function.

(ii) $\mu(\Omega)=0$.}

\medskip

Equivalently, one can define analytic pseudo--measures as finitely
additive functions on the Boolean algebra generated by a basis 
of clopen sets for the topology of $\Omega$ satisfying the conditions
of Definition 5.1.1.

\smallskip

In particular, we will be interested in the case where the 
space $\Omega=\partial \cT$ is the boundary of a tree. 
One defines currents on trees in the following way.

\medskip

{\bf 5.1.2. Definition.} {\it 
Let $W$ be an abelian group and let $\cT$ be a locally finite tree. 
We denote by $\cC(\cT,W)$ the
group of $W$--valued {\it currents} on $\cT$. These are $W$--valued maps
$\bold{c}$ from the set of oriented edges of $\cT$ to $W$ that satisfy 
the following properties:

\smallskip

(i) Orientation reversal: $\bold{c}(\bar e)=-\bold{c}(e)$, where $\bar e$ 
denotes the edge $e$ with the reverse orientation.

\smallskip

(ii)  Momentum conservation: 
$$
\sum_{s(e)=v} \bold{c}(e)= 0, \eqno(5.1)
$$
where $s(e)$ (resp. $t(e)$) denote the source (resp. target) vertex of the
oriented edge $e$.}

\medskip

One can then identify currents on a tree with analytic pseudo--measures 
on its boundary, as in [vdP]. 

\medskip

{\bf 5.1.3. Lemma.} {\it
The group $\cC(\cT,W)$ of currents on $\cT$ is canonically isomorphic
to the group
of $W$--valued finitely additive analytic pseudo--measures on $\partial \cT$.}

\medskip

{\bf Proof.} The identification is obtained by setting
$$
\mu(V(e))= \bold{c}(e), \eqno(5.2)
$$
where $V(e)\subset \partial\cT$ is the clopen subset of the boundary
of $\cT$ defined by all infinite admissible (i.e.~ without backtracking)
paths in $\cT$ starting with the oriented edge $e$.

\medskip

The group of $W$--valued currents on $\cT$ can also be characterized in
the following way. We let $\cA(\cT,W)$ denote the $W$--valued functions
on the oriented edges of $\cT$ satisfying $\mu(\bar e)=-\mu(e)$ and
let $C(\cT^{(0)},W)$ denote the $W$--valued functions on the set of
vertices of $\cT$. The following result is also proved in [vdP].

\medskip

{\bf 5.1.4. Lemma.} {\it
Let $d: \cA(\cT,W)\to C(\cT^{(0)},W)$ be given by
$$ d(f)(v)=\sum_{s(e)=v} f(e). $$
Then the group of $W$--valued currents is given by
$\cC(\cT,W)=Ker(d)$, so that one has an exact sequence
$$ 0 \to \cC(\cT,W)\to \cA(\cT,W) \to C(\cT^{(0)},W) \to
0. $$}

\medskip

{\bf 5.2. The tree of $PSL(2,\Z)$ and its boundary.}
We will apply these results to the tree $\cT$ 
of $PSL(2,\Z)$ embedded in the 
hyperbolic plane $\H$. Its vertices are the elliptic points 
$\widetilde I\cup \widetilde R$, where $\widetilde I$ is  the 
$PSL(2,\Z)$--orbit of  $i$ 
and $\widetilde R$ is the orbit of $\rho=e^{2\pi i/3}$.
The set of edges is given by the geodesic arcs
$\{ \gamma (i),\gamma (\rho) \}$, for $\gamma \in SL(2,\Z)$.

\smallskip

The relation between $\partial \cT$ and $\P^1(\R)$
can be summarized as follows. 

\medskip

{\bf 5.2.1. Lemma.} {\it
The boundary $\partial \cT$ is a compact Hausdorff space. There is a
natural continuous $PSL(2,\Z)$--equivariant surjection 
$\Upsilon: \partial \cT \to \partial \H=\P^1(\R)$, 
which is one--to--one on the irrational points $\P^1(\R)\cap
(\R\smallsetminus \Q)$ and two--to--one on rational points.}

\medskip

{\bf Proof.}
Consider the Farey tessellation of the hyperbolic plane $\H$ by
$PSL(2,\Z)$ translates of the ideal triangle with vertices $\{
0,1,\infty\}$. The tree $\cT\subset \H$ has a vertex of valence three
in each triangle, a vertex of valence two bisecting each edge of the
triangulation and three edges in each triangle joining the valence
three vertex to each of the valence two vertices on the sides of the
triangle.

\smallskip 

If we fix a base vertex in the tree, for instance the vertex
$v=\rho=e^{2\pi i/3}$, then we can identify the boundary $\partial
\cT$ with the set of infinite admissible paths (i.e.~ paths without
backtracking) in the tree $\cT$ starting at $v$.

\smallskip 

Any such path
traverses an infinite number of triangles and can be encoded by the
sequence of elements $\gamma_n\in PSL(2,\Z)$ that determine the
successive three--valent vertices crossed by the path, $v_n= \gamma_n
v$. 

\smallskip

This sequence of points $v_n\in \H$ accumulates at some point
$\theta \in \P^1(\R)=\partial \H$. If the point is irrational, 
$\theta\in \P^1(\R)\smallsetminus \P^1(\Q)$, then the point $\theta$
is not a vertex of any of the Farey triangles and there is a unique
admissible sequence of vertices $v_n\in \tilde R\subset \H$ with the
property that $\lim_{n\to \infty} v_n =\theta$. To see this, 
consider the family of segments $I_n\subset \P^1(\R)$ given by
all points that are ends of an admissible path in $\cT$ starting at
$v_n$ and not containing $v_{n-1}$.

\smallskip 

The intersection $\cap_n I_n$
consists of a single point which is identified with the point in the
limit set $\P^1(\R)$ of $PSL(2,\Z)$ given by the infinite sequence 
$\gamma_1\gamma_2\cdots \gamma_n \cdots$. 

\input epsf
\midinsert
$$\centerline{\hbox{\epsffile{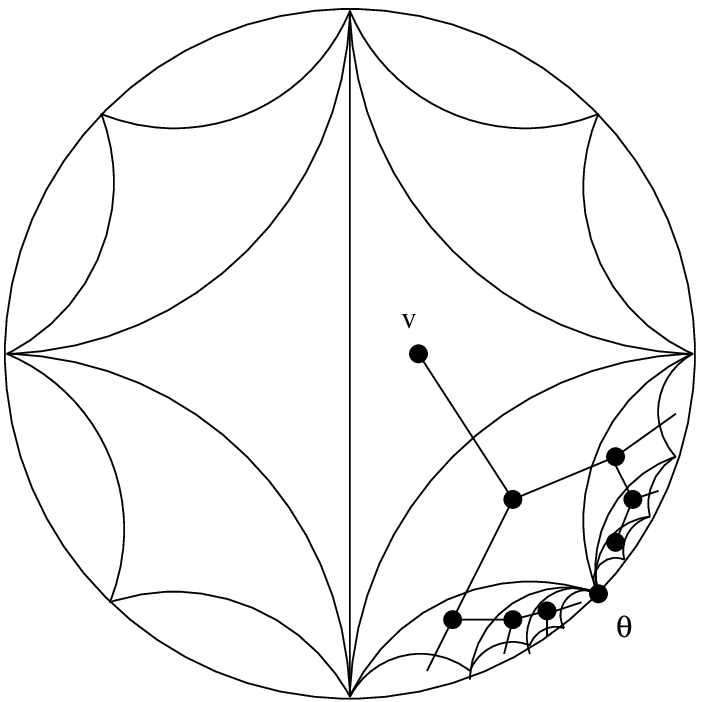}}}$$
\endinsert

Consider now a rational point $\theta\in \P^1(\Q)$. Then
$\theta$ is a vertex of some (in fact infinitely many) of the 
Farey triangles. In this case, one can see that there are two distinct
sequences of $3$-valent vertices $v_n$ of $\cT$ with the property that
$\lim_{n\to \infty} v_n =\theta$, due to the fact that, beginning with
the first $v_k$ such that $\theta$ is a vertex of the triangle
containing $v_k$, there are two adjacent triangles that also have 
$\theta$ as a vertex, see Figure 1.

\smallskip

This defines a map $\Upsilon: \partial \cT \to \partial \H$, given by
$\Upsilon(\{ v_n \})=\lim_{n\to \infty} v_n =\theta$. By construction
it is continuous, $1:1$ on the irrationals and $2:1$ on the
rationals. 

\medskip

{\bf 5.3. The disconnection space.} 
An equivalent description of the space $\partial \cT$ of ends of the
tree $\cT$ of $PSL(2,\Z)$ can be given in terms of the {\it disconnection}
spaces considered in [Spi]. We discuss it here briefly, as it
will be useful later in describing the noncommutative geometry of the
boundary action of $PSL(2,\Z)$.

\smallskip

Given a subset $U\subset\P^1(\R)$ one considers the abelian
$C^*$-algebra $\cA_U$ generated by the algebra $C(\P^1(\R))$ and the
characteristic functions of the positively oriented intervals 
(in the sense of \S 1.1 above) with endpoints in $U$. If the set $U$
is dense in $\P^1(\R)$ then this is the same as the closure in the
supremum norm of the $*$-algebra generated by these
characteristic functions.

\smallskip

The Gelfand--Naimark correspondence $X\leftrightarrow C_0(X)$
furnishes an equivalence of categories
of locally compact Hausdorff topological spaces $X$ and commutative
$C^*$--algebras respectively. Thus, we have
$\cA_U=C(D_U)$, where the topological space $D_U$ is called {\it the
disconnection} of $\P^1(\R)$ along $U$. It is a compact Hausdorff space
and it is totally disconnected if and only if $U$ is dense in
$\P^1(\R)$.

\smallskip

In particular, one can consider the set $U=\P^1(\Q)$ and the
resulting disconnection $D_\Q$ of $\P^1(\R)$ along $\P^1(\Q)$. 

\medskip

{\bf 5.3.1. Lemma.}\, {\it
The map $\Upsilon: \partial \cT \to \P^1(\R)$ of Lemma 5.2.1
factors through a homeomorphism $\widetilde\Upsilon:  \partial \cT \to D_\Q$, 
followed by the surjective map $D_\Q \to
\P^1(\R)$ determined by the inclusion of the algebra $C(\P^1(\R))$
inside $C(D_\Q)$. }

\medskip

{\bf Proof.} 
The compact Hausdorff space $\partial \cT$ is dual to the
abelian $C^*$-algebra $C(\partial\cT)$. After the choice of a base
vertex $v\in \cT$, the topology of $\partial \cT$ is generated 
by the clopen sets $V(v')$, of ends of admissible paths starting at a
vertex $v'$ and not passing through $v$. 

\smallskip

In fact, it suffices to consider only
$3$--valent vertices, because for a vertex $v'$ of valence two we have
$V(v')= V(v'')$ with $v''$ being the next $3$--valent vertex in the direction
away from $v$. The $C^*$--algebra $C(\partial\cT)$ is thus
generated by the characteristic functions of the $V(v')$, since
$\partial \cT$ is a totally disconnected space. Since 
$v'=\gamma v$ for some $\gamma\in PSL(2,\Z)$, the set
$V(v')=\Upsilon^{-1}(I(v')) \subset \partial \cT$, where $I(v')\subset
\P^1(\R)$ is a segment of the form $I(v')=\gamma I$, where $I$ is
one of the segments $[\infty,0]$ or $[0,1]$ or $[1,\infty]$. These
are the segments in
$\P^1(\R)$ of the form $[p/q,r/s]$ for $p,q,r,s\in \Z$ with $ps-qr=\pm
1$, that is, the primitive segments on the boundary, as defined in \S 1.6,
corresponding to sides of triangles of the Farey tessellation. 

\smallskip

On the other hand the disconnection $D_\Q$ of $\P^1(\R)$ along
$\P^1(\Q)$ is dual to the abelian $C^*$-algebra generated by 
all the characteristic functions of the oriented intervals 
with endpoints in $\alpha, \beta\in\P^1(\Q)$. 
Since the sets $V(e)$ that generate the topology of $\partial\cT$
are of this form, this shows that there is an injection 
$C(\partial\cT) \to C(D_\Q)$. The primitive intervals of \S 1.6 
in fact give a basis for the topology of $D_\Q$, since the characteristic
functions that generate $C(D_\Q)$ can be written as combinations
of characteristic functions of such intervals, using primitive
chains as in \S 1.6. This shows that the map is in fact an isomorphism. 
 
\medskip

{\bf 5.4. Analytic pseudo--measures on $D_\Q$.} 
We show here that the pseudo--measures on $\P^1(\R)$ 
(Definition 2.1.1) can be regarded as 
analytic pseudo--measures on the disconnection space $D_\Q$. 

\medskip

{\bf 5.4.1. Lemma.}\, {\it
There is a natural bijection between $W$--valued pseudo--measures on $\P^1(\R)$ and $W$--valued analytic pseudo--measures on $D_\Q$ 
(or equivalently, $W$--valued
currents on $\cT$.) }

\medskip

{\bf Proof.}
Recall that a basis for the topology of $D_\Q$ is given 
by the preimages under the map $\Upsilon$ of the primitive (Farey) segments 
in $\P^1(\R)$. These are segments of the form $(g(\infty),g(0))$ 
for some $g\in G=PSL(2,\Z)$. Let $e$ be an oriented edge in the tree
$\cT$ of $PSL(2,\Z)$. A basis of the topology of $\partial\cT$ is
given by sets of the form $V(e)$. These are in fact Farey intervals
and all such intervals arise in this way. Thus, an analytic
pseudo--measure on $D_\Q$ is a finitely additive function on the
Boolean algebra generated by the Farey intervals with the properties
of Definition 5.1.1.

\smallskip

A pseudo--measure on $\P^1(\R)$, on the other hand, is  a map
$\mu: \P^1(\Q)\times \P^1(\Q) \to W$ with the properties (1.1) and (1.2).
Such a function in fact extends, as we have seen, to 
a finitely additive function on the Boolean algebra consisting of finite
unions of (positively oriented) intervals $(\alpha,\beta)$ with
$\alpha,\beta\in \P^1(\Q)$.

\smallskip

We then just set $\mu_{an}(V(e)):=\mu(\alpha,\beta)$, where $\P^1(\R)\supset
[\alpha,\beta] =\Upsilon (V(e))$. This is then a finitely additive
function on the 
Boolean algebra of the $V(e)$. The condition
$\mu(\beta,\alpha)=-\mu(\alpha,\beta)$ 
implies that $\mu_{an}(\partial\cT)=0$ and similarly this combined
with the property 
$\mu(\alpha,\beta)+\mu(\beta,\gamma)+\mu(\gamma,\alpha)=0$ implies that
$\mu_{an}(V\cup V')=\mu_{an}(V)+\mu_{an}(V')$ for $V\cap V' =\emptyset$.

\smallskip

Conversely, start with an analytic pseudo--measure $\mu_{an}$ on $D_\Q$
with values in $W$. For $\Upsilon V(e)=[\alpha,\beta]$ we write 
$\mu(\alpha,\beta):=\mu_{an}(V(e))$. Then the properties 
$\mu_{an}(V\cup V')=\mu_{an}(V)+\mu_{an}(V')$ for $V\cap V' =\emptyset$ and
$\mu_{an}(\partial\cT)=0$ imply the corresponding properties
$\mu(\beta,\alpha)=-\mu(\alpha,\beta)$ 
and $\mu(\alpha,\beta)+\mu(\beta,\gamma)+\mu(\gamma,\alpha)=0$ for the
pseudo--measure. 

\smallskip

Under this correspondence, $\Gamma$--modular $W$--valued pseudo--measures
correspond 
to $\Gamma$--modular analytic $W$--valued pseudo--measures. Thus, in the
following we often 
use the term pseudo--measure equivalently for the analytic ones on
$D_\Q$ or for those defined 
in \S 2 on $\P^1(\R)$.

\medskip

{\bf 5.5. K--theoretic interpretation.}
In [MaMar1] we gave an interpretation of the modular complex of
[Ma1] in terms of $K$-theory of $C^*$-algebras, by considering
the crossed product $C^*$-algebra of the action of $G=PSL(2,\Z)$ on 
its limit set $\P^1(\R)$. 
Here we consider instead the action of $PSL(2,\Z)$ on the ends
$\partial\cT=D_\Q$ of the tree and we obtain a similar $K$-theoretic
interpretation of pseudo--measures. 

\smallskip

More generally, let $\Gamma \subset PSL(2,\Z)$ be a finite index
subgroup and we consider the action of $G=PSL(2,\Z)$ on 
$\partial\cT\times \P$, where $\P=\Gamma\backslash
PSL(2,\Z)$. We let $\cA=C(\partial\cT\times \P)$ with the 
action of $G=PSL(2,\Z)=\Z/2\Z* \Z/3\Z$ by automorphisms. 

\medskip

{\bf 5.5.1. Lemma.} {\it The $K$--groups of
the crossed product $C^*$--algebra $\cA\rtimes G$
have the following structure:
$$
K_0(\cA\rtimes G)=Coker(\alpha), \ \ \  K_1(\cA\rtimes
G)=Ker(\alpha), 
$$
where 
$$ 
\alpha: C(\partial\cT\times \P,\Z)\to C(\partial\cT\times
\P,\Z)^{G_2}\oplus C(\partial\cT\times
\P,\Z)^{G_3} 
$$
is given by $\alpha: f \mapsto (f+f\circ \sigma, f+f\circ \tau+f
\circ \tau^2)$ for $\sigma$ and $\tau$, respectively, the generators
of $G_2$ and $G_3$ defined by (2.4).}

\medskip

{\bf Proof.} 
First recall that, for a totally disconnected space $\Omega$, 
one can identify the locally constant integer valued 
functions $C(\Omega,\Z)$ with $K_0(C(\Omega))$, 
whereas $K_1(C(\Omega))=0$.

\smallskip 

The six--terms exact sequence of [Pim] for groups
acting on trees gives
$$ 0 \to K_1(\cA\rtimes G) \to  K_0(\cA) \to K_0(\cA\rtimes
G_2)\oplus K_0(\cA\rtimes G_3) \to K_0(\cA\rtimes G) \to 0, $$
for $G =PSL(2,\Z)$ and $G_i$ equal to $\Z/2\Z$ or $\Z/3\Z$.
Let $\alpha: K_0(\cA\rtimes
G_2)\oplus K_0(\cA\rtimes G_3) \to K_0(\cA\rtimes G)$ be the
map in the sequence above.

\smallskip

Here we use the fact that $K_1(\cA)=0$ and $K_1(\cA\rtimes
G_i)=K^1_{G_i}(\cA)=0$ so that the remaining terms in the six--terms
exact sequence do not contribute.

\smallskip

Moreover, we have $K_0(\cA)=C(\partial\cT\times \P,\Z)$. 
Similarly, we have 
$$ K_0(\cA\rtimes G_i)=K^0_{G_i}(\partial\cT\times
\P)=C(\partial\cT\times \P,\Z)^{G_i} . $$ 
Following \cite{Pim}, we see that the map $\alpha$ can be described as
the map 
$$ 
\alpha: f \mapsto (f+f\circ \sigma , f+f\circ \tau+f
\circ \tau^2). 
$$

\smallskip

This gives a description of $Ker(\alpha)$ as 
$$ 
Ker(\alpha) =\{ f\in C(\partial\cT\times \P,\Z)\,|\, 
f+f\circ \sigma = f+f\circ \tau+f \circ \tau^2=0 \}. 
$$
Similarly, the cokernel of $\alpha$ is the group of coinvariants, that
is, the quotient of $C(\partial\cT\times \P,\Z)$ by the submodule generated
by the elements of the form $f+f\circ \sigma$ and $f+f\circ \tau+f
\circ \tau^2$. By the six--terms exact sequence we know that
$K_1(\cA\rtimes G)=Ker(\alpha)$ and $K_0(\cA\rtimes G)=
Coker(\alpha)$. 

\smallskip

For simplicity let us reduce to the case with $\P=\{1 \}$, that is,
$\Gamma=G=PSL(2,\Z)$. 

\medskip

{\bf 5.5.2. Integral.} 
Given a $W$--valued pseudo--measure $\mu$ on $D_\Q$, in the sense
of Definition 5.1.1,  we can  define an integral 
$$ 
f\mapsto \int f \, d\mu \, \in W, 
$$
for $f\in C(D_\Q,\Z)$ in the following way.

\smallskip

An element $f\in C(D_\Q,\Z)$ is of the form
$f=\sum_{i=1}^n a_i \chi_{I_i}$ with $a_i\in \Z$ and the
intervals $I_i\subset D_\Q$ with $\Upsilon I_i= g_i (\infty,0) \subset
\P^1(\R)$, for some $g_i\in G$. Thus, a natural prescription is 
$$ 
\int f \, d\mu =\sum_i a_i \mu(I_i) \in W. 
$$ 

\smallskip

To simplify notation, in the following 
we often do not distinguish between an interval in 
$\P^1(\R)$ and its lift to the disconnection space $D_\Q$. 
So we write equivalently $\mu(I_i)$ or $\mu(g_i(\infty),g_i(0))$. 

\medskip

{\bf 5.5.3. Lemma.} {\it Let $g\in G$. 
The following  change of variable
formula holds:
$$ 
\int f\circ g\, d\mu = \int f\, d(\mu\circ g^{-1}). 
$$}

\smallskip

{\bf Proof.} 
We have $f\circ g =\sum_i a_i \chi_{I_i}\circ g=
\sum_i a_i \chi_{g^{-1}I_i}$, so that
$$ \int f\circ g\, d\mu = \sum_i a_i
\mu(g^{-1}g_i(\infty),g^{-1}g_i(0)) $$
$$ = \sum_i a_i \mu\circ g^{-1} (g_i(\infty),g_i(0))= \int f\, d(\mu\circ
g^{-1}). $$

\smallskip

The following generalization is also true (and easy):
$$ 
\int (f\circ g)\, h\, d\mu = \int f\,(h\circ g^{-1}) 
d(\mu\circ g^{-1}), 
$$
for $f,h\in C(D_\Q,\Z)$. 

\medskip

{\bf 5.5.4. Proposition.} {\it Let $\mu$ be
 a $G$-modular $W$--valued pseudo--measure.
 For any $h\in K_1(\cA\rtimes G)$, there exists 
 a unique $G$-modular $W$-valued pseudomeasure $\mu_h$ with
$$ 
\mu_h(\infty,0) := \int h\, d\mu  
$$
}

\medskip

{\bf Proof.} We will consider the case of $G=PSL(2,\Z)$.
The general case can be reduced to this one by
proceeding as in  2.4 for the modular pseudomeasures.

\smallskip

In view of the Theorem 2.3, it suffices to check that
the element
$$
\int h\, d\mu  \in W, 
$$
is annihilated by $1+\sigma$ and
$1+\tau+\tau^2$.

\smallskip

An element $h\in K_1(\cA\rtimes G)$ is a
function $f\in C(\partial\cT,\Z)$ satisfying
$h+h\circ\sigma=0$ and $h+h\circ\tau+h\circ\tau^2=0$.
Therefore
$$ 
(1+\sigma) \int h\, d\mu =\int h\,d\mu + \int h\,d\mu\circ\sigma
= \int (h+h\circ\sigma) d\mu =0. 
$$
Similarly,
$$ 
(1+\tau+\tau^2) \int h\, d\mu = \int (h+h\circ\tau^2+h\circ\tau)
d\mu=0. 
$$
Thus, one obtains a $G$-modular pseudo-measure $\mu_h$ for each
$h\in K_1(\cA\rtimes G)$.

\medskip

{\bf 5.5.5. Proposition.} {\it A $G$-modular $W$--valued pseudo--measure $\mu$
defines a group homomorphism $\mu: K_0(\cA\rtimes G) \to W$
induced by integration with respect to $\mu$.}

\medskip

{\bf Proof.} Using
the identification with analytic pseudomeasures, we can
consider the functional from $C(\partial\cT,\Z)$ to $W$ given by
integration
$$ 
\mu(f)=\int f \, d\mu . 
$$
This descends to the quotient of $C(\partial\cT,\Z)$ by the relations
$f+f\circ \sigma$ and $f+f\circ\tau+f\circ \tau^2$. In fact,
it suffices to consider the case where $f$ is the characteristic
function of a segment $(g(\infty),g(0))$. We have
$$
\int (f+f\circ\sigma) d\mu= \int f\,d\mu + \int f \,d\mu\circ \sigma
=(1+\sigma)\int f\,d\mu  
$$
and
$$ \int \chi_{(g(\infty),g(0))} \, d\mu = \int \chi_{(\infty,0)}
\circ g^{-1} \, d\mu = \int \chi_{(\infty,0)} \, d\mu\circ g $$
so that
$$ \sigma \int \chi_{(\infty,0)} \, d\mu\circ g =
\int \chi_{(\infty,0)} d\mu\circ g \circ \sigma =
\int \chi_{(g\sigma(\infty),g\sigma(0))} \, d\mu, $$
hence
$$ (1+\sigma)\int \chi_{(g(\infty),g(0))} \, d\mu =
g(1+\sigma)\mu(\infty,0) =0. $$
The argument for $f+f\circ\tau+f\circ\tau^2$ is similar.

\medskip

{\bf 5.6. Boundary action and noncommutative spaces.}
In [MaMar1], we described the noncommutative boundary of the
modular tower in terms of the quotient 
$\Gamma\backslash \P^1(\R)$, which we interpreted as a noncommutative
space, described by the crossed product $C^*$-algebra 
$C(\P^1(\R))\rtimes \Gamma$ or $C(\P^1(\R)\times \P)\rtimes G$, for
$\Gamma\subset G$ a finite index subgroup and $\P$ the coset space.

\smallskip

Here we have seen that, instead of considering $\P^1(\R)$ as the
boundary, one can also work with the disconnection space $D_\Q$.
We have then considered the corresponding crossed product
$C(D_\Q\times \P)\rtimes G$. This is similar to the treatment of
Fuchsian groups described in [Spi].

\smallskip

There is, however, another way of describing the boundary action of
$G=PSL(2,\Z)$ on $\P^1(\R)$. As we have also discussed at length
in [MaMar1], it uses a dynamical system associated to the
Gauss shift of the continued fraction expansion, generalized in the
case of $\Gamma\subset G$ to include the action on the coset space
$\P$. In [MaMar1] we worked with $PGL(2,\Z)$ instead of
$PSL(2,\Z)$. For the $PSL(2,\Z)$ formulation, see [ChMay], 
[Mayer] and [KeS]. 

\smallskip

This approach to describing the boundary geometry was at the basis of
our extension of modular symbols to limiting modular symbols. We
discuss briefly here how this formulation also leads to a
noncommutative space, in the form of an Exel--Laca algebra.

\medskip

{\bf 5.6.1. The generalized Gauss shift.}
The Gauss shift for $PSL(2,\Z)$ is the map $\hat T:[-1,1]\to [-1,1]$
of the form
$$
 \widehat T: x\mapsto -sign(x) T(|x|),  
$$
where
$$ 
T(x)=\frac{1}{x} - \left[ \frac{1}{x}\right].
$$
Notice that this differs from the Gauss shift $T:[0,1]\to
[0,1]$, $T(x)=1/x -[1/x]$ of $PGL(2,\Z)$ by the presence of the extra sign,
as in [ChMay], [KeS].

\smallskip

When one considers a finite index subgroup $\Gamma\subset G=PSL(2,\Z)$
one extends the shift $\hat T$ to the generalized Gauss shift
$$
\hat T_\P : \cI \times \P \to \cI \times \P, \ \ \ (x,s)\mapsto (\hat
T(x), [g ST^k]),
$$
where here $\cI=[-1,1]\cap (\R\smallsetminus \Q)$ and where
$k=sign(x)n_1$. Here $\P=\Gamma\backslash G$ is the coset space and $g\in G$
denotes the representative $\Gamma g=s\in \P$.
The $n$--th iterate of the map $\hat T_\P$ acts on $\P$
as the $SL(2,\Z)$ matrix
$$ \left( \matrix -sign(x_1) p_{k-1}(x) & (-1)^k p_k(x) \\
q_{k-1}(x) & (-1)^{k+1} q_k(x) \endmatrix \right). $$

\medskip

{\bf 5.6.2. Shift happens.}
Instead of considering the action of the generalized Gauss shift on the space $\cI\times\P$, one 
can proceed as in [KeS] and introduce a shift space over which $\hat T_\P$ acts as a
shift operator. This is obtained by considering the countable alphabet $\Z^\times\times \P$ and
the set of {\it admissible sequences} 
$$
\Sigma_\Gamma=\{ ((x_1,s_1),(x_2,s_2),\ldots)\, |\,
A_{(x_i,s_i),(x_{i+1},s_{i+1})}=1 \},
$$
where the matrix $A$ giving the admissibility condition is defined as
follows. One has $A_{(x,s),(x',s')}=1$ if $xx'<0$ and
$$ s'=\tau_x(s):= [g ST^{x_1}] \in \P, \ \ \text{ where } s=\Gamma g $$
and $A_{(x,s),(x',s')}=0$ otherwise. The
action of $\widehat T_\P$ described above becomes the action of the
one--sided shift $\sigma: \Sigma_\Gamma \to \Sigma_\Gamma$
$$
 \sigma: ((x_1,s_1),(x_2,s_2),\ldots) \mapsto
((x_2,s_2),(x_3,s_3),\ldots) .
$$

The space $\Sigma_\Gamma$ can be topologized in the way that is customarily
used to treat this type of dynamical systems given by shift
spaces. Namely, one considers on $\Sigma_\Gamma$ the topology generated by
the cylinders (all words in $\Sigma_\Gamma$ starting with an assigned
finite admissible word in the alphabet). This makes $\Sigma_\Gamma$ into a
compact Hausdorff space. One can see, in fact, that the topology is
metrizable and induced for instance by the metric
$$
d((x_k,s_k)_k, (x'_k,s'_k)_k)=\sum_{n=1}^\infty 2^{-n}
\left(1-\delta_{(x_n,s_n), (x'_n,s'_n)}\right).
$$

As is shown in [KeS], the use of this topology as opposed to the
one induced by $\P^1(\R)\times \P$ simplifies the analysis of the
associated Perron--Frobenius operator. The latter now falls into the
general framework developed in [MaUr] and one obtains the
existence of a shift invariant ergodic measure on $\Sigma_\Gamma$ from this
general formalism. One uses essentially the {\it finite
irreducibility} of the shift space $(\Sigma_\Gamma,\sigma)$, which follows
in [KeS] from the ergodicity of the geodesic flow.
The $\sigma$--invariant measure on $\Sigma_\Gamma$ induces via the bijection
between these spaces a $\hat T_\P$--invariant measure on the space $\cI
\times \P$.

\medskip

{\bf 5.6.3. Exel--Laca algebras.}
There is a way to associate to a shift space on an alphabet a 
noncommutative space, in the form of a Cuntz--Krieger algebra
in the case of a finite alphabet [CuKrie], or more 
generally an Exel--Laca algebra for a countable alphabet [ExLa].

\smallskip 

Start with a (finite or countable) alphabet
$\fA$ and a matrix $A: \fA\times \fA \to \{ 0,1 \}$ that assigns
the admissibility condition for words in the alphabet $\fA$.
In the case of a finite alphabet $\fA$, the corresponding
Cuntz--Krieger algebra is the $C^*$-algebra
generated by partial isometries $S_a$ for $a\in \fA$ with the
relations
$$
S_a S_a^* = P_a, \ \ \ \text{ with } \sum_a P_a =1   \eqno(5.3)
$$
$$
S_a^* S_a = \sum_{b\in\fA} A_{ab} S_b S_b^*.      \eqno(5.4)
$$
In the case of a countably infinite alphabet $\fA$, one has to be more
careful, as the summations that appear in the relations (5.3)
and (5.4) no longer converge in norm. A version of CK
algebras for infinite matrices was developed by Exel and Laca in
[ExLa]. One modifies the relations (5.3)
and (5.4) in the following way.

\medskip

{\bf 5.6.4. Definition ([ExLa]).} {\it
For a countably infinite alphabet $\fA$ and an admissibility matrix
$A:\fA\times \fA \to \{ 0,1 \}$, the CK algebra $O_A$ is the
universal $C^*$-algebra generated by partial isometries $S_a$, for
$a\in \fA$, with the following conditions:
\smallskip

(i) $S_a^* S_a$ and $S_b^* S_b$ commute for all $a,b\in \fA$.
\smallskip
(ii) $S_a^* S_b=0$ for $a\neq b$.
\smallskip
(iii) $(S_a^* S_a) S_b = A_{ab} S_b$ for all $a,b\in \fA$.
\smallskip
(iv) For any pair of finite subsets $X,Y\subset \fA$, such that the
product
$$
A(X,Y,b):= \prod_{x\in X} A_{xb} \prod_{y\in Y} (1-A_{yb}) 
\eqno(5.5)
$$
vanishes for all but finitely many $b\in \fA$, one has the identity
$$
\prod_{x\in X} S_x^* S_x \prod_{y\in Y} (1-S_y^* S_y) =
\sum_{b\in \fA} A(X,Y,b) S_b S_b^* .  
\eqno(5.6)
$$}

\smallskip

The conditions listed above are obtained by formal manipulations from
the relations (5.3) and (5.4) and are equivalent to
them in the finite case.

\medskip

{\bf 5.6.5. The algebra of the generalized Gauss shift.}
We now consider the shift space $(\Sigma_\Gamma,\sigma)$ of \cite{KeStr}.
In this case we have the alphabet $\fA=\Z^\times \times \P$ and the
admissibility matrix given by the condition
$$
A_{ab}=1 \ \ \ \text{ iff } nn'<0 \ \ \text{ and } \ \  s=s'\, ST^{n'},
$$
for $a=(n,s)$ and $b=(n',s')$ in $\fA$. The corresponding Exel--Laca algebra
$O_A$ is generated by isometries $S_{(n,s)}$ satisfying the conditions of
Definition 5.5.4 above.

\smallskip

Let $\cJ=\Upsilon^{-1}[-1,1] \subset D_\Q$ be the preimage of the
interval $[-1,1]$ under the continuous surjection $\Upsilon:D_\Q
\to \P^1(\R)$ and let $\cJ_{+1} =\Upsilon^{-1}[0,1]$ and
$\cJ_{-1}=\Upsilon^{-1}[-1,0]$. Also let
$$
\cJ_k=\{ x\in \cJ \,|\, \Upsilon(x)=sign(k) [a_1,a_2,\ldots], \,
a_1=|k|, \, a_i\geq 1\}.
$$

Let $\Gamma$ be a finite index subgroup of $G=PSL(2,\Z)$,
$\P=\Gamma\backslash G$. Consider the sets
$$
\cJ_{k,s}:= \cJ_k \times \{s\},
$$
for $s\in \P$ and for $k\in \Z^\times$. Let $\chi_{k,s}$
denote the characteristic function of the set $\cJ_{k,s}$.

\medskip

{\bf 5.6.6. Proposition.} {\it
The subalgebra of the crossed product $C(D_\Q\times \P)\rtimes G$
generated by elements of the form
$$
S_{k,s}:= \chi_{k,s} U_k, 
$$
where $U_k=U_{\gamma_k}$, for 
$\gamma_k=T^k S \in \Gamma$,
is isomorphic to the Exel--Laca algebra $O_A$ of the shift
$(\Sigma_\Gamma,\sigma)$.}

\medskip

{\bf Proof.}
We need to check that the $S_{k,s}$ satisfy the Exel--Laca
axioms of Definition 5.5.4. We have
$S_{k,s}^* = U_k^* \chi_{k,s}= U_k^*(\chi_{k,s}) U_k^*$
so that $S_{k,s} S_{k,s}^*= \chi_{k,s} U_k U_k^* \chi_{k,s}
=\chi_{k,s}=: P_{k,s}$ and $S_{k,s}^* S_{k,s}=U_k^*
\chi_{k,s}\chi_{k,s} U_k= U_k^*(\chi_{k,s})$. One sees from this that
$S_{k,s}^*S_{k,s}$ and $S_{k',s'}^*S_{k',s'}$ commute, as they both
belong to the subalgebra $C(D_{\P^1(\Q)}\times \P)$, so that the first
axiom is satisfied. Similarly, $S_{k,s}^*S_{k',s'}= U_k^* \chi_{k,s}
\chi_{k',s'} U_{k'}=0$ for $(k,s)\neq (k',s')$, which is the second
axiom. To check the third axiom, we have
$$ (S_{k,s}^*S_{k,s})  S_{k',s'} = U_k^*(\chi_{k,s}) \chi_{k',s'}
U_{k'}. 
$$

We now describe more explicitly the element $U_k^*(\chi_{k,s})$.
Since we have $U_k^*(f)=f\circ \gamma_k$, we consider the
action of the element $\gamma_k= T^k S\in
PSL_2(\Z)$ on the set $\cJ_k \times \{ s \}$.
We have $S(\cJ_k)=\{ x\,|\, \pi(x)= -k -sign(k) [a_2,a_3,\ldots] \}$,
and $T^k S(\cJ_k)=\{ x\,|\, \pi(x)=-sign(k) [a_2,a_3,\ldots]
\}$. This is the union of all $$\{ x \,|\,
\pi(x)=sign(k')[|k'|,a_3,\ldots]\,\, sign(k')=-sign(k) \}.$$
Thus, we have
$$
T^k S: \cJ_k \times \{ s \} \to \cup_{k': kk'<0}\,\, \cJ_{k'} \times \{ s\,
ST^{-k} \}, \eqno(5.7)
$$
we obtain that
$$
U_k^*(\chi_{k,s}) \chi_{k',s'} =A_{(k,s),(k',s')} \chi_{k',s'}. \eqno(5.8)
$$
Thus, we obtain $(S_{k,s}^*S_{k,s})  S_{k',s'} = A_{(k,s),(k',s')}
S_{k',s'}$, which is the third Exel--Laca axiom.

\smallskip

For the last axiom, consider the condition that
$A(X,Y,b)$ of (5.5) vanishes for all but finitely many
$b\in \fA$. Given two finite subsets $X,Y\subset \Z^\times \times \P$,
the only way that $A(X,Y,b)=0$ for all but finitely many
$b\in \Z^\times \times \P$ is that $A(X,Y,b)=0$ for all $b$.

\smallskip

In fact, for given $X$ and $Y$, suppose that there exists an element
$b=(k,s)$ such that $A(X,Y,b)\neq 0$. This means that, for all $x\in
X$ we have $A_{xb}=1$ and for all $y\in Y$ we have $A_{yb}=0$. The
first condition means that for all $x=(k_x,s_x)$ we have $k_x k<0$ and
$s=s_x\, ST^{k_x}$, while the second condition means that, for all
$y=(k_y,s_y)$ we either have $k_y k>0$ or $k_y k <0$ and $s\neq s_y\,
ST^{k_y}$. Consider then elements of the form $b'=(k',s)$ with
$k'k>0$ and the same $s\in \P$ as $b=(k,s)$. All of these still
satisfy $A_{xb'}=1$ for all $x\in X$ and $A_{yb}=0$ for all $y\in Y$,
since the conditions only depend on the sign of $k$ and on $s\in
\P$. Thus, there are infinitely many $b'$ such that $A(X,Y,b')\neq 0$.
If $A(X,Y,b)\equiv 0$ for all $b$, the condition of the fourth axiom
reduces to
$$
 \prod_{x\in X} S_x^* S_x \prod_{y\in Y} (1-S_y^*S_y) =0. \eqno(5.9)
$$

Suppose that $A(X,Y,b)\equiv 0$ for all $b\in \fA$ but that
the expression in (5.9) is non-zero. This means
that both $\prod_{x\in X} S_x^* S_x\neq 0$ and $\prod_{y\in Y}
(1-S_y^*S_y)\neq 0$. We first show by induction that, for $X=\{ x_1,
\ldots, x_N\}$ we can write first product in the form
$$
(\prod_{x\in X} S_x^* S_x) P_b = A_{x_1 b} A_{x_2,b} \cdots A_{x_N b} P_b, \eqno(5.10)
$$
where $P_b =\chi_b =S_b S_b^*$. We know already that this is true for
a single point $X=\{ x \}$ by (5.8). Suppose it is true for
$N$ points. Then we have
$$ S_{x_0}^* S_{x_0} (\prod_{i=1}^N S_{x_i}^* S_{x_i})  P_b =
(\prod_{i=1}^N A_{x_i b})  \,\, S_{x_0}^* S_{x_0} \,\,  P_b =
A_{x_1 b} A_{x_2,b} \cdots A_{x_N b} \,\, A_{x_0 b} P_b, $$
which gives the result.
Thus, we see that the condition $\prod_{x\in X} S_x^* S_x\neq 0$
implies that, for some $b\in \fA$, one has $A_{x,b}\neq 0$ for all
$x \in X$, i.e.~ one has $\prod_{x\in X} A_{x,b}\neq 0$. We analyze
similarly the condition $\prod_{y\in Y} (1-S_y^*S_y)\neq 0$. We show
by induction that, for $Y=\{ y_1, \ldots, y_M \}$, the product can be
written in the form
$$
(\prod_{y\in Y} (1-S_y^*S_y)) P_b = (1- A_{y_1 b}) (1- A_{y_2 b})\cdots
(1-A_{y_M b}) P_b. 
\eqno(5.11)
$$
First consider the case of a single point $Y=\{ y \}$. We have
$1- S_y^*S_y = 1- U_{k_y}^*(\chi_{k_y,s_y})$ so that using
(5.8) we get $(1- S_y^*S_y) P_b = (1-A_{yb}) P_b$.
We then suppose that the identity (5.11) holds for $M$
points. We obtain, again using (5.8),
$$
(1-S_{y_0}^*S_{y_0}) (\prod_{i=1}^M (1- A_{y_i b}))  P_b =
(\prod_{i=1}^M (1- A_{y_i b}))\,\, (1-A_{y_0 b}) P_b.
$$
Thus, the condition $\prod_{y\in Y}(1- S_y^*S_y)\neq 0$ implies that
there exists $b\in \fA$ such that $\prod_{y\in Y} (1-A_{yb})\neq 0$,
which contradicts the fact that we are assuming $A(X,Y,b)\equiv 0$.
This proves that the fourth Exel--Laca axiom is satisfied.

\medskip

To summarize, there are several interesting noncommutative spaces 
related to the boundary of modular curves $X_\Gamma$: the crossed
product algebras $C(\P^1(\R)\times \P)\rtimes G$ and $C(D_\Q\times
\P)\rtimes G$, and the Exel--Laca algebra $O_A$. This calls for a more 
detailed investigation of their relations and of the information 
about the modular tower that each of these noncommutative spaces 
captures.

\medskip

{\bf 5.7. Limiting modular pseudo--measures.}
In [MaMar1] we extended classical modular symbols to include the case 
of limiting cycles associated to geodesics with irrational endpoints, the limiting
modular symbols. The theory of limiting modular symbols was further studied
in [Mar] and more recently in [KeS].

\smallskip

We show here that we can  similarly define limiting
modular pseudo--measures, extending the class of pseudo--measures from 
finitely additive functions on the Boolean algebra of the $V(e)\subset D_\Q$ to
a larger class of sets by a limiting procedure.

\smallskip

In the following we assume that the group $W$ where the pseudo--measures take values
is also a real (or complex) vector space. If as a vector space $W$ is infinite dimensional,
we assume that it is a topological vector space. This is always the case for finite dimension.
As before we also assume that $W$ is a left $G$-module.

\smallskip

We recall the following property of the convergents of the continued fraction expansion, cf.~ 
[PoWe]. For $\theta\in \R\smallsetminus \Q$, let $q_n$ denote, as before, the successive 
denominators of the continued fraction expansion. Then the limit
$$ 
\lambda(\theta):= \lim_{n\to \infty} \frac{ 2 \log q_n(\theta) }{n} 
\eqno(5.12) 
$$
is defined away from an exceptional set $\Omega\subset \P^1(\R)$ and where it
exists it is equal to the Lyapunov exponent of the Gauss shift. 

\smallskip

In [MaMar1] we proved that the limiting modular symbol 
$$
\{\{*,\theta\}\}= \lim_{t\to\infty}\frac{1}{t} \{ x_0, y(t) \}_\Gamma
$$
 can be computed by the limit
$$  
\lim_{n\to\infty} \frac{1}{\lambda(\theta)n} \sum_{k=1}^{n+1} \{ g_{k-1}(0),g_{k-1}(\infty)\}_\Gamma
 \in H_1(X_\Gamma,\R) 
\eqno(5.13) 
$$
where $g_k=g_k(\theta)$ is the matrix in $G$ that implements the action of the 
$k$--th power of the Gauss shift.

\smallskip 

Fix a $W$--valued pseudo--measure $\mu$ on $\P^1(\R)$, equivalently thought
of as an analytic pseudo--measure on $D_\Q$.  Consider the class of 
positively oriented 
intervals $(\infty,\theta)$ with $\theta\in \P^1(\R)\smallsetminus \P^1(\Q)$ 
and define the limit
$$ 
\mu^{lim}(\infty,\theta):=  \lim_{n\to\infty} \frac{1}{\lambda(\theta) n} \sum_{k=1}^{n+1} 
\mu(g_{k-1}(\infty),g_{k-1}(0)) \in W.  
\eqno(5.14) 
$$
This is defined away from the exceptional set $\Omega$ in $\P^1(\R)$ which contains $\P^1(\Q)$ as
well as the irrational points where either the limit defining $\lambda(\theta)$ does not exist
or the limit of the $$ \frac{1}{n}\sum_{k=1}^{n+1} \mu(g_{k-1}(\infty),g_{k-1}(0)) $$ does not exist.

\smallskip 

Similarly, for $\theta\in \P^1(\R)\smallsetminus \Omega$ we set
$$ 
\mu^{lim}(\theta,\infty)= \lim_{n\to\infty} \frac{1}{\lambda(\theta) n} \sum_{k=1}^{n+1} 
\mu(g_{k-1}(0,\infty)). 
\eqno(5.15) 
$$
We then set
$$ \mu^{lim}(\theta,\eta):= \mu^{lim}(\theta,\infty)+\mu^{lim}(\infty,\eta). 
\eqno(5.16) 
$$

\medskip

{\bf 5.7.1. Lemma.} {\it
The function $\mu^{lim}$ defined as above satisfies
$$ 
\mu^{lim}(\eta,\theta)=-\mu^{lim}(\theta,\eta) \ \ \ \ 
\mu^{lim}(\theta,\eta)+\mu^{lim}(\theta,\zeta)+\mu^{lim}(\zeta,\theta)=0 
$$
for all $\theta,\eta,\zeta\in \P^1(\R)\smallsetminus \Omega$.}

\medskip

{\bf Proof.}
Since $(\beta,\alpha)=g(0,\infty)=g\sigma(\infty,0)$, we have
$\mu(\beta_k,\alpha_k)=\mu(g_k(0,\infty))=\mu(g_k\sigma(\infty,0))=-\mu(\alpha_k,\beta_k)$ so that
$$ \mu^{lim}(\theta,\infty)=-\mu^{lim}(\infty,\theta). $$
We then have $\mu^{lim}(\eta,\theta)=\mu^{lim}(\infty,\theta)+\mu^{lim}(\eta,\infty)=
-\mu^{lim}(\theta,\eta)$. The argument for the second identity is similar.

\medskip

Thus, the limiting pseudo--measure $\mu^{lim}$ defines a finitely additive function on the
Boolean algebra generated by the intervals $(\theta,\eta)$ with $\theta,\eta\in \P^1(\R)\smallsetminus \Omega$. The limiting pseudo--measure $\mu^{lim}$ is $G$--modular if $\mu$ is $G$--modular.

\medskip

{\bf 5.7.2. Limiting pseudo--measures and currents.}
In terms of currents on the tree $\cT$, we can describe the limit computing $\mu^{lim}$ as 
a process of averaging the current $\bold{c}$ over edges along a path. 

\medskip

{\bf 5.7.3. Lemma.} {\it
Let $\mu$ be a $W$--valued pseudo--measure and $\bold{c}$ the corresponding current on $\cT$.
The limiting pseudo--measure $\mu^{lim}$ is computed by the limit
$$ \mu^{lim}(\infty,\theta) = \lim_{n\to \infty}\frac{1}{\lambda(\theta) n} \sum_{k=1}^n \bold{c}(e_k). $$}

\medskip

{\bf Proof.}
An irrational 
point $\theta$ in $\partial\cT$ corresponds to a unique admissible
infinite path in the tree $\cT$  
starting from a chosen base vertex. We describe such a path as an
infinite admissible sequence 
of oriented edges $e_1 e_2\cdots e_n \cdots$. To each such edge there
corresponds an open set  
$V(e_k)$ in $\partial\cT$ with the property that $\cap_k V(e_k)=\{
\theta \}$. These correspond  
under the map $\Upsilon: \partial\cT \to \P^1(\R)$ to intervals
$[g_k(\infty),g_k(0)]$. Thus, the  
expression computing the limiting pseudo--measure can be written
equivalently via the limit of the  averages
$$ \frac{1}{n} \sum_{k=1}^n \bold{c}(e_k). $$

\bigskip

\centerline{\bf References}

\medskip

[A] A.~Ash. {\it Parabolic cohomology of arithmetic subgroups
of $SL(2,\Z )$ with coefficients in the field of rational functions
on the Riemann sphere.} Amer. Journ. of Math., 111 (1989), 35--51.

\smallskip

[ChMay] C.H.~Chang, D.~Mayer. {\it Thermodynamic formalism and
Selberg's zeta function for modular groups}. Regular and Chaotic
Dynamics 15 (2000)  N.3, 281--312.

\smallskip

[ChZ] Y.~J.~Choie, D.~Zagier. {\it Rational period functions
for for $PSL(2,\Z ).$} In: A Tribute to Emil Grosswald:
Number Theory and relatied Analysis, Cont. Math.,
143 (1993), AMS, Providence, 89--108.

\smallskip

[CuKrie] J.~Cuntz, W.~Krieger. {\it A class of $C^*$--algebras
and topological Markov chains, I,II}. Invent. Math. 56 (1980)
251--268 and Invent. Math. 63 (1981) 25--40.

\smallskip

[ExLa] R.~Exel, M.~Laca. {\it Cuntz--Krieger algebras for
infinite matrices}.  J. Reine Angew. Math.  512  (1999), 119--172.

\smallskip

[Fu1] Sh.~Fukuhara. {\it Modular forms, generalized Dedekind symbols
and period polynomials.} Math. Ann., 310 (1998), 83--101.

\smallskip

[Fu2] Sh.~Fukuhara. {\it Hecke operators on weighted Dedekind symbols.}
Preprint math.NT/0412090.

\smallskip

[He1] A.~Herremans. {\it A combinatorial interpretation of 
Serre's conjecture on modular Galois representations.}
PhD thesis, Catholic University Leuven, 2001.

\smallskip

[He2] A.~Herremans. {\it A combinatorial interpretation of 
Serre's conjecture on modular Galois representations.}
Ann. Inst. Fourier, Grenoble, 53:5 (2003), 1287--1321.

\smallskip

[KeS] M.~Kessenb\"ohmer, B.~O.~Stratmann. {\it Limiting modular symbols
and their fractal geometry.} Preprint math.GT/0611048.

\smallskip

[Kn1] M.~Knopp. {\it Rational period functions of the modular group.}
Duke Math. Journal, 45 (1978), 47--62.

\smallskip

[Kn2] M.~Knopp. {\it Rational period functions of the modular group II.}
Glasgow Math. Journal, 22 (1981), 185--197.

\smallskip

[L] P.~L\'evy. {\it Sur les lois de probabilit\'e dont
d\'ependent les quotients complets et incomplets d'une fraction
continue.} Bull.~Soc.~Math.~France, 557 (1929), 178--194.

\smallskip

[LewZa]  J.~Lewis, D.~Zagier. {\it Period functions and
the Selberg zeta function for the modular group.}
In: The Mathematical Beauty of Physics, Adv. Series in
Math. Physics 24, World Scientific, Singapore, 1997, pp. 83--97. 

\smallskip

[Ma1] Yu.~Manin. {\it Parabolic points and zeta-functions of modular curves.}
Russian: Izv. AN SSSR, ser. mat. 36:1 (1972), 19--66. English:
Math. USSR Izvestiya, publ. by AMS, vol. 6, No. 1 (1972), 19--64,
and Selected Papers, World Scientific, 1996, 202--247.

\smallskip

[Ma2] Yu.~Manin. {\it Periods of parabolic forms and $p$--adic Hecke series.}
Russian: Mat. Sbornik, 92:3 (1973), 378--401. English:
Math. USSR Sbornik, 21:3 (1973), 371--393,
and Selected Papers, World Scientific, 1996, 268--290.

\smallskip

[Ma3] Yu.~Manin. {\it Iterated integrals of modular forms and 
noncommutative modular symbols.}
In: Algebraic Geometry and Number Theory.
In honor of V. Drinfeld's 50th birthday.
Ed. V.~Ginzburg. Progress in Math., vol. 253. 
Birkh\"auser, Boston, pp.565--597. Preprint math.NT/0502576.

\smallskip

[Ma4] Yu.~Manin. {\it Iterated Shimura integrals.} Moscow Math. Journal, 
vol. 5, Nr. 4 (2005), 869--881.
 Preprint math.AG/0507438.

\smallskip

[MaMar1] Yu.~Manin, M.~Marcolli. {\it Continued fractions, modular symbols, and non-commutative geometry.} Selecta math., new ser. 8 (2002),
475--521. Preprint math.NT/0102006.

\smallskip

[MaMar2] Yu.~Manin, M.~Marcolli. {\it Holography principle and arithmetic of algebraic curves.} Adv. Theor. Math. Phys., 5 (2001), 617--650.
Preprint hep-th/0201036.

\smallskip

[Mar] M.~Marcolli. {\it Limiting modular symbols and the
Lyapunov spectrum}. Journal of Number Theory, 98 (2003) 348--376.

\smallskip

[MaUr] R.D.~Mauldin, M.~Urba\'nski. {\it Graph directed Markov
systems. Geometry and dynamics of limit sets}. Cambridge Tracts in
Mathematics, 148. Cambridge University Press, Cambridge, 2003.

\smallskip

[Mayer] D.~Mayer. {\it Continued fractions and related
transformations}.  Ergodic theory, symbolic dynamics, and hyperbolic
spaces (Trieste, 1989),  175--222, Oxford Sci. Publ., Oxford Univ.
Press, New York, 1991.

\smallskip

[Me] L. Merel. {\it Universal Fourier expansions of modular forms.} In:
Springer LN in Math., vol. 1585 (1994), 59--94.

\smallskip

[Pim] M.V.~Pimsner, {\it $KK$--groups of crossed products by
groups acting on trees}.  Invent. Math.  86  (1986),  no. 3, 603--634.  

\smallskip

[PoWe] M.~Pollicott, H.~Weiss. {\it Multifractal analysis of 
Lyapunov exponent for continued fraction and Manneville-Pomeau 
transformations and applications to Diophantine approximation}.  
Comm. Math. Phys.  207  (1999),  no. 1, 145--171. 

\smallskip

[Sh1] V.~Shokurov. {\it The study of the homology of Kuga
varieties.} Math. USSR Izvestiya, 16:2 (1981), 399--418.

\smallskip

[Sh2] V.~Shokurov. {\it Shimura integrals of cusp forms.}
Math. USSR Izvestiya, 16:3 (1981), 603--646.

\smallskip

[Spi] J.S.~Spielberg. {\it Cuntz--Krieger algebras associated
with Fuchsian groups}. Ergodic Theory Dynam. Systems 13 (1993), no.
3, 581--595.

\smallskip

[vdP] M.~van der Put. {\it Discrete groups, Mumford curves and
theta functions}. Ann. Fac. Sci. Toulouse Math. 6 (1992) 399--438.

\enddocument